\documentclass[12pt,reqno]{amsart}

\usepackage{amssymb}
\usepackage{amsmath}
\usepackage{bbm}
\usepackage{mathrsfs}

\makeatletter
\renewcommand{\@biblabel}[1]{#1.}
\makeatother

\theoremstyle{plain}
\newtheorem{prp}{Proposition}
\newtheorem{lemma}[prp]{Lemma}
\newtheorem{thm}[prp]{Theorem}
\newtheorem{cor}[prp]{Corollary}

\theoremstyle{remark}
\newtheorem{remark}{Remark}
\numberwithin{prp}{section}

\newcommand{\rf}[1][]{\textup{\eqref{#1}}}

\newcommand\ov{\overline}
\newcommand\wbar{{\overline{w}}}
\newcommand{\wo}{{w_\circ}}
\newcommand{\wokk}{{w_\circ^\prime}}

\newcommand{\Rda}[2][]{\hat{R}_{#1}^{#2}}
\newcommand{\Rch}[2][]{\check{R}_{#1}^{#2}}
\newcommand{\Rre}[1][]{\acute{R}_{#1}}
\newcommand{\Rli}[2][]{\grave{R}_{#1}^{#2}}
\newcommand{\Rdam}[1][]{\hat{R}^{-1}_{#1}}

\newcommand{\Rchm}[1][]{\check{R}^{-1}_{#1}}

\newcommand{\Rrepm}[1][]{\acute{R}^\tau_{#1}}

\DeclareMathOperator{\id}{id}
\DeclareMathOperator{\End}{End}
\DeclareMathOperator{\Mor}{Mor}
\DeclareMathOperator{\im}{im}
\DeclareMathOperator{\adr}{Ad_\rr}
\DeclareMathOperator{\rank}{rank}
\DeclareMathOperator{\codim}{codim}

\DeclareMathOperator{\tr}{tr}
\DeclareMathOperator{\TR}{Tr}
\DeclareMathOperator{\Card}{Card}

\newcommand\dd{\mathrm{d}}

\newcommand\ot{\otimes}
\newcommand{\ott}{{\otimes}}
\newcommand{\ota}{{\otimes}_{\!{\scriptscriptstyle\mathcal{A}}}}

\newcommand{\ww}{{\mbox{$\scriptscriptstyle{W}$}}}

\newcommand\wt{\widetilde}
\newcommand\til{\wt{T}_w}
\newcommand\Rel{\mathrm{Rel}}
\newcommand{\fr}[2][]{\mathfrak{#2}_{#1}}

\newcommand\AAAA{\boldsymbol{A}}
\newcommand\AAA{\boldsymbol{a}}
\newcommand\bsig{{\boldsymbol\sigma}}
\renewcommand\SS{\mathscr{S}}
\newcommand{\ant}[1][k]{A_{#1}}

\newcommand{\tz}{{\tau,z}}
\newcommand\ntz{\mathfrak{n}_\tz}
\renewcommand{\S}[1][k]{\mathcal{S}_{#1}}
\renewcommand\AA{\mathcal{A}}
\newcommand\B{\mathcal{B}}
\newcommand\CC{\mathcal{C}}
\newcommand{\dett}{\mathcal{T}}
\newcommand\DD{\mathcal{D}}
\newcommand\PP{\mathcal{P}}
\newcommand\RR{\mathcal{R}}

\renewcommand\H{{\mathcal H}_k(q)}

\newcommand{\h}[1][k]{H_{#1}(q)}
\newcommand\hh{H_k}

\newcommand\C{\mathbbm{C}}
\newcommand\N{\mathbbm{N}}
\newcommand\NO{\mathbbm{N}_0}

\newcommand{\slqn}[1][N]{SL_q(#1)}
\newcommand{\glqn}[1][N]{GL_q(#1)}
\newcommand{\oslqn}[1][N]{\mathcal{O}(SL_q(#1))}
\newcommand{\oglqn}[1][N]{\mathcal{O}(GL_q(#1))}

\newcommand\nn{\nonumber}
\newcommand\ds{\displaystyle} 
\newcommand\ts{\textstyle}

\newcommand{\qm}[1][1]{q^{-{#1}}}
\newcommand{\QM}{\hat{q}}
\newcommand{\QP}{\check{q}}

\newcommand{\lam}{{\lambda}}
\newcommand\vp{\varphi}
\newcommand\eps{\epsilon}
\newcommand\ve{\varepsilon}
\newcommand\vt{\vartheta}
\newcommand\vr{\varrho}
\newcommand\vN{\varTheta}
\newcommand\vO{\varOmega}
\newcommand\om{\omega}

\newcommand\rac{{\,\triangleleft\,}}
\newcommand{\fettc}{{\scriptstyle{\tt c}}}
\newcommand{\fettt}{{\scriptstyle{\tt t}}}
\newcommand{\inv}  {{\scriptstyle{\rm inv}}}

\newcommand{\uhr}{{\upharpoonright}}

\newcommand\ii{\mathrm{i}}
\renewcommand\ll{\ell}
\newcommand\rr{\mathrm{r}}
\newcommand\dl{{\Delta_\ll}}
\newcommand\dr{{\Delta_\rr}}

\newcommand\half{{\frac{1}{2}}}
\newcommand\ivec{{\vec{\imath}}}    
\newcommand\jvec{{\vec{\jmath}}}   
\newcommand\nvec{{\vec{n}}}        
\newcommand\xvec{{\vec{x}}}
\newcommand\xinv{{\overset{\leftarrow}{x}}}
\newcommand\zvec{{\vec{z}}}

\newcommand\iinv{{\overset{\leftarrow}{\imath}}}  
\newcommand\jinv{{\overset{\leftarrow}{\jmath}}}

\newcommand{\Gamm}{\varGamma}
\newcommand\gpmz{\Gamm_{\pm,z}}
\newcommand\gtz{\Gamm_{\tau,z}}

\newcommand\gl{{\Gamm_\ll}}
\newcommand\gr{{\Gamm_\rr}}
\newcommand\gi{{\Gamm_\ii}}
\newcommand{\gd}[1][]{\Gamm^{\land #1}} 
\newcommand\sig{\sigma}
\newcommand\sigg{\gamma}
\newcommand{\gten}[1][]{\Gamm^{\otimes #1}}
\newcommand{\gdw}[1][]{{_\ww{\Gamm^{\land #1}}}}
\newcommand{\gds}{{_s{\Gamm^{\land}}}}
\newcommand{\gdu}{{_u{\Gamm^{\land}}}}
\newcommand\ul{\underline}
\newcommand\JJ{{_\ww J}}
\newcommand\Js{{_s J}}
\newcommand\Ju{{_u J}}
\newcommand\lo{{\ell(\wo)}}

\newcommand{\gdsl}[1][]{{_s{\Gamm^{\land #1}_\ll}}}

\newcommand{\gdwi}[1][]{{_\ww{\Gamm^{\land {#1}}_\ii}}}
\newcommand{\gdwl}[1][]{{_\ww{\Gamm^{\land {#1}}_\ll}}}
\newcommand{\gol}[1][]{{\Gamm^{\ot {#1}}_\ll}}
\newcommand\pair{{\langle\cdot,\cdot\rangle}}
\newcommand\uckuk{(u^\fettc)^{\ot k}\ott u^{\ot k}}
\newcommand\ucuk{(u^\fettc\ott u)^{\ot k}}
\newcommand{\ip}[2][]{I_{#1}^{(#2)}}
\newcommand{\idp}[1]{I^{\land #1}}
\newcommand{\ipm}[2][]{I_{\tau\, {#1}}^{(#2)}}

\begin{document}

\author{Axel Sch\" uler}

\thanks{Supported by the Deutsche Forschungsgemeinschaft, e-mail: schueler@mathematik.uni-leipzig.de}

\address{Axel Sch\" uler, Department of Mathematics, 
University of Leipzig, Augustusplatz 10, 04109 Leipzig,
Germany}

\title[Differential Hopf Algebras]{Differential Hopf Algebras \\
on Quantum
  Groups of Type A}

\begin{abstract}
Let $\AA$ be a Hopf algebra and $\Gamm$ be a bicovariant first order
differential calculus over $\AA$. It is known that there
are three possibilities to construct
a differential  Hopf algebra $\gd=\gten /J$  that  contains $\Gamm$ as its
first order part. 
Corresponding to the three choices of the ideal $J$,
we distinguish the `universal' exterior algebra, 
the `second antisymmetrizer' exterior algebra, and Woronowicz' external
algebra, respectively.
\\
Let $\Gamm$ be one
of the $N^2$-dimensional bicovariant first order differential calculi on
the quantum group   $\glqn$ or $\slqn$, and let $q$ be a transcendental complex
number. For  Woronowicz' external algebra we determine the dimension
of the space of  left-invariant and
of bi-invariant $k$-forms, respectively.   Bi-invariant forms are  closed and
represent different  de Rham cohomology classes. The algebra of bi-invariant
forms is graded anti-commutative.
\\
For $N\ge3$ the three differential Hopf algebras coincide.
However, in case of the $4D_\pm$-calculi on  $SL_q(2)$
the universal differential Hopf algebra is strictly larger than Woronowicz'
external algebra. The bi-invariant 1-form is not closed.
\end{abstract}

\maketitle


\section{Introduction}

Non-commutative differential geometry on quantum groups is 
a basic tool for further applications in both theoretical 
physics and mathematics. A general framework for bicovariant differential
calculus on quantum groups has been invented by Woronowicz \cite{Wo2}.
Covariant  first order differential calculi (abbreviated FODC) were
constructed, studied, and classified by many authors, see (for instance)
\cite{CSWW,J,Su,SS2,HSm}. Despite the rather extensive 
literature on bicovariant {\em first}\/ order differential calculi
the corresponding exterior algebras
have been treated only in few cases, see \cite{Scha,Mal}.
The de Rham cohomology of the three dimensional left-covariant differential
calculus  on the quantum group $SU_q(2)$  was calculated by Woronowicz
\cite{Wo1}. The  de Rham cohomology  of the four dimensional bicovariant
differential calculi $4D_\pm$ on  $SU_q(2)$ were calculated by Grie\ss l
\cite{G}. Brzezi\'nski \cite{B} pointed out that the exterior algebra $\gds$
is a differential Hopf algebra.
\\
The purpose of this paper is to compare three  possible constructions 
of differential Hopf algebras (exterior algebras) over quantum groups of type~A. Let $\AA$ be a Hopf algebra
and let $\Gamm$ be a bicovariant FODC over $\AA$. Consider the tensor
algebra $\gten$ over $\AA$. Let $\Ju$ be  the two-sided ideal generated by
the elements $\sum_{(r)}\om(r_{(1)})
\ota\om(r_{(2)})$, $r\in \RR$, where 
$\om(a)=\sum_{(a)}Sa_{(1)}\dd a_{(2)}$ and
$\RR=\ker\ve\cap\ker\om$ is the associated right ideal.
Let $\Js$ denote  the ideal generated by $\ker(I-\sig)$, where $\sig$ is the braiding  of $\Gamm\ota\Gamm$.
Finally let  $\JJ=\bigoplus_{k\ge 2}\ker \ant$, where $\ant$ 
is the $k$th antisymmetrizer
constructed from the braiding $\sig$.
Define the exterior algebras $\gdu=\gten/\Ju$, $\gds=\gten/\Js$,
 and $\gdw=\gten/\JJ$
and call them {\em universal exterior algebra},
{\em second antisymmetrizer exterior algebra}\/   and {\em Woronowicz' 
external algebra}, respectively.
The first one is the ``largest" one. It can be characterized
by the  following universal property: Each  differential Hopf algebra
with a given FODC $\Gamm$ as its first order part is a quotient of $\gdu$, see
\cite[Subsect.\,14.3.3]{KS} or \cite[Theorem\,5.5]{LS}.
Both the second and the third constructions  use the braiding $\sig$.
This is a  twisted flip automorphism
of the bicovariant bimodule $\Gamm\ota\Gamm$.  It
satisfies the braid equation.  The second antisymmetrizer
exterior algebra uses the antisymmetrizer $\ant[2]=I-\sig$ only. 
The definition  of Woronowicz' external algebra however involves antisymmetrizers of all degrees.
\\
Now let $\AA$ be one of the Hopf algebras $\oglqn$ or
$\oslqn$. Let $\Gamm=\gtz$, $\tau\in\{+,-\}$, $z\in\C$, $z\ne0$,
 denote one of the $N^2$-dimensional bicovariant FODC over 
$\AA$ constructed in \cite{SS3} by a method
of Jur\v co. Our standing assumptions are that $q$  is a transcendental complex
number and $N\ge2$. We present three main results. The first two are exclusively
concerned with Woronowicz' construction while the last one compares the three
possible exterior algebras. 
The first  result is stated in Theorem\,\ref{left}. It says that the dimension
of the space of left-invariant $k$-forms equals $\binom{N^2}{k}$. In particular, there is  a unique up to scalars  left invariant form of maximal degree $N^2$
(this form is even bi-invariant). The second main result stated in 
Theorem\,\ref{bi} is concerned with the subalgebra of bi-invariant forms. This
algebra is graded anti-commutative.
The dimension of the space of bi-invariant $k$-forms is equal to the number of
partitions of $k$ into a sum of pairwise different  positive odd integers less
than $2N$. Bi-invariant forms are closed and represent different  de Rham cohomology classes.
The third main result is stated in Theorem\,\ref{usw}.  The differential Hopf
algebras  $\gds$ and $\gdw$ are isomorphic. Suppose the
 parameter $z$  be regular (only finitely many values of $z$ are excluded).
 Then $\gdu$ and $\gds$ are isomorphic differential Hopf  algebras. For
 $\AA=\oglqn[2]$ or $\AA=\oslqn[2]$, $\Gamm=\Gamm_{+,z}$, and  $z^2=q^{-2}$
 however, the universal differential calculus $\gdu$ is strictly larger than $\gds$. The bi-invariant
 1-form $\theta\in\gdu$ is not closed and $\theta^2$ is central.
\\
The paper is organised as follows. In Section\,\ref{pre} we recall 
preliminary facts  about bicovariant bimodules,  bicovariant
first order and  higher order differential calculus over Hopf algebras;
we also recall the  construction of bicovariant FODC on $\glqn$ and  $\slqn$.  In
Section\,\ref{result} we formulate the main results. Section\,\ref{s-hecke} is
devoted to the Iwahori-Hecke algebra $\h$. Our first two main results  are obtained by studying the `abstract' $\sig$-algebra. In case
of quantum groups of type~A  the braiding $\sig$ can be identified with
$T^{-1}\ott T$, where $T$ is a single generator of $\h$.  In this way Woronowicz' antisymmetrizer can be viewed as an element of
$\h\ott\h$. The basic result, given in
Proposition\,\ref{spectrum}, describes the  decomposition of the `abstract'
antisymmetrizer as a linear combination of mutually orthogonal idempotents
$\pi_\lam$. In Section\,\ref{s-corep} we recall facts from the theory of
corepresentations of  $\AA=\oglqn$ and $\AA=\oslqn$. We essentially make use
of our assumption that $q$ is transcendental. Firstly, $\AA$ is then
cosemisimple and the tensor product of corepresentations decomposes into
irreducible ones. Secondly, this decomposition is as in the classical case. 
By Brauer-Schur-Weyl duality idempotents of  $\h$ correspond to
subcorepresentations of $u^{\ot k}$ where $u$ denotes the fundamental matrix
corepresentation of $\AA$.  The basic result is Proposition\,\ref{p-dec}. It
says that $\pi_\lam$ induces a corepresentation of $\AA$ that is isomorphic
to the tensor product of two single irreducible ones. The
connection between the `abstract' antisymmetrizer and  Woronowicz'
antisymmetrizer is recovered in Proposition\,\ref{p-mor}; Theorem\,\ref{left}
and Theorem\,\ref{bi} are proved here. The proof of Theorem\,\ref{usw} is given in
Section\,\ref{s-univers}. The first part of the proof uses facts about dual
quadratic algebras, see \cite{Ma}, and transmutation theory invented by Majid,
see \cite{Maj}. The second and third parts of the proof are very technical
since we need the explicit description of the associated right ideal.
\section{Preliminaries}\label{pre}
Throughout the paper we work over the ground field $\C$ (with one  exception
in the proof of Lemma\,\ref{biinv-l}
where the field $\C(q)$ of rational functions is used). All vector spaces,
algebras, bialgebras, etc.\ are meant to be $\C$-vector spaces, unital $\C$-algebras,
$\C$-bialgebras etc. The linear span of a set $\{a_i\colon i\in K\}$ is
denoted by $\langle a_i\colon i\in K\rangle$.
In this paper $\AA$ denotes a bialgebra or a Hopf algebra.
All modules, comodules, and bimodules are assumed to be $\AA$-modules, 
$\AA$-comodules, and $\AA$-bimodules if nothing else is specified. 
We denote the comultiplication, the counit, and the
antipode by $\Delta$, $\ve$, and  by $S$, respectively.
We use the notions ``right comodule" and ``corepresentation" of $\AA$ as 
synonyms. By fixing a basis in the underlying vector space we identify
corepresentations and the corresponding matrices.  
Let $v$ (resp.\ $f$) be a corepresentation (resp.\  a representation) of
$\AA$. As usual $v^\fettc$ (resp.\  $f^\fettc$)  
denotes the contragredient corepresentation (resp.\  contragredient
representation) of $v$ (resp.\  of $f$). The
space of intertwiners of corepresentations $v$ and $w$ is  $\Mor(v,w)$.
We write $\Mor(v)$
for $\Mor(v,v)$. By $\End$ and $\ot$ we always mean
$\End_\C$ and $\ot_\C$, respectively. If $A$ is a linear mapping, $A^\fettt$
denotes the transpose of $A$ and $\tr A$ the trace of $A$.
Lower indices of $A$ always refer
to the components of a tensor product where $A$ acts (`leg numbering').
The unit matrix is
denoted by $I$.
  We set $\wt{a}=a-\ve(a)$ for
$a\in\AA$.
We use Sweedler's notation for the coproduct $\Delta(a)=\sum a_{(1)}
\ot a_{(2)}$,
for left comodules  $\vp(e)=\sum e_{(-1)}\ot e_{(0)}$, and for right
comodules  $\psi(e)=\sum e_{(0)}\ot e_{(1)}$.
The mapping $\adr\colon\AA\to\AA\ot\AA$,
 $\adr a=\sum a_{(2)}\ot Sa_{(1)} a_{(3)}$,
is a right comodule map called the {\em right adjoint coaction} 
\/of $\AA$ on itself.
\subsection*{ Bicovariant bimodules and tensor algebra}
A {\em bicovariant bimodule}\/ over $\AA$ (or  {\em Hopf bimodule}\/)
is a bimodule $\Gamm$ together with linear mappings 
$\dl\colon\Gamm\to\AA\ot\Gamm$ and $\dr\colon\Gamm\to\Gamm\ot\AA$ such that
$(\Gamm,\dl)$ is a left comodule, $(\Gamm,\dr)$ is a right comodule,
$(\id\ot\dr)\dl=(\dl\ot\id)\dr$,
$\dl(a\omega b)=\Delta(a)\dl(\omega)\Delta(b)$, 
and $\dr(a\omega b)=\Delta(a)\dr(\omega)\Delta(b)$
for $a,\,b\in\AA$ and $\omega\in\Gamm$.
 An element
$\omega\in\Gamm$ is called {\em left-invariant}\/ (resp.\  {\em right-invariant}\/)
if $\dl(\omega)=1\ot\omega$ (resp.\  $\dr(\omega)=\omega\ot 1$).
The linear space
of left-invariant (resp.\  right-invariant) elements of $\Gamm$ is denoted by $\gl$ 
(resp.\  $\gr$). The elements of $\gi=\gl\cap\gr$ are called {\em bi-invariant}.
The structure of bicovariant bimodules has been  completely characterized 
by Theorems 2.3 and 2.4 in \cite{Wo2}. We recall the corresponding result:
\\
Let $(\Gamm,\dl,\dr)$ be a bicovariant bimodule over $\AA$ and let 
$\{\omega_i:i\in K\}$ be a finite linear basis of $\gl$. 
Then there exist matrices $v=(v^i_j)$
and $f=(f^i_j)$ of elements $v^i_j\in \AA$ and of functionals $f^i_j$ on
$\AA$, $i,j\in K$ such that: 
\\
(i) $\sum_{(a)} Sa_{(1)}\omega_i a_{(2)}
=\sum_n f^i_n(a)\omega_n$, $a\in \AA$, $i\in K$,
 and $\dr(\omega_i)=\sum_n\omega_n\ot v^n_i$, $i\in K$.
\\
(ii) $v$ is a corepresentation and $f$ is a representation of $\AA$. 
\\
(iii) $\sum_r v^r_i(a\ast f^r_j)=
\sum_r(f^i_r\ast a)v^j_r$, $a\in\AA$, $i,j\in K$.
 \\
We have set $a\ast f= \sum f(a_{(1)})a_{(2)}$ and $f\ast a=\sum
a_{(1)}f(a_{(2)})$. The set $\{\omega_i:i\in K\}$ is a free left module basis of $\Gamm$.
Conversely, if $\{\omega_i: i\in K\}$
is a basis of a certain finite dimensional vector space 
$\Gamm_0$ and if $v$ and $f$  are matrices satisfying (ii)
and (iii) then there exists a unique bicovariant bimodule $\Gamm$ such 
that $\gl=\Gamm_0$ and (i) holds. In this situation we simply 
write $\Gamm=(v,f)$. 
\\
Let  $\Gamm=(v,f)$, $\Gamm_1=(v_1,f_1)$, and $\Gamm_2=(v_2,f_2)$ be bicovariant bimodules.
It is easy to check that the tensor product of $\AA$-bimodules 
$\Gamm_1\ota\Gamm_2$ is also a  bicovariant bimodule with bimodule structure, 
left coaction $\dl$, and right coaction $\dr$ defined by
 $a{\cdot}\alpha\ota\beta {\cdot} b=
a\alpha\ota\beta b$,
$\dr(\alpha\ota\beta)=\sum\alpha_{(0)}\ota\beta_{(0)}\ot
\alpha_{(1)}\beta_{(1)}$, and
$\dl(\alpha\ota\beta)=\sum\alpha_{(-1)}\beta_{(-1)}
\ot\alpha_{(0)}\ota\beta_{(0)}$.
The corresponding pair of corepresentation and
representation is $(v_1\ott v_2, f_1\ott f_2)$. Similarly,
$\gten[k]=\Gamm\ota\cdots\ota\Gamm$ ($k$ factors), $k\ge2$, and
$\gten=\bigoplus_{k\ge 0}\gten[k]$, where
$\Gamm^{\ot 0}=\AA$ (by definition, $\dl=\dr=\Delta$ on $\AA$),  $\Gamm^{\ot 1}=\Gamm$, are  bicovariant bimodules.
Note that  $\Gamm_\ll$ becomes  a right module via $\rho\rac a=\sum Sa_{(1)} \rho a_{(2)}$, $a\in \AA$. Moreover
$(\gten)_\ll$ is right module algebra.
\\
{\em The shuffle decomposition.} 
\/Let $\S$ be the symmetric group on $\{1,\dots,k\}$ and let
$s_i=(i,i+1)\in\S$, $i=1,\dots, k-1$,
denote the simple transposition that exchanges $i$ and $i+1$. 
Each $w\ne 1$ in $\S$ can be written 
in the form $w=s_{i_1}\cdots s_{i_r}$ for some $i_n\in\{1,\dots,k-1\}$.
If $r$ is as small as
possible call it the {\em length}\/ of $w$, written $\ell(w)$, and call
any expression of $w$ as the product of $r$ elements of
 $\{s_1,\dots,s_{k-1}\}$ a 
{\em reduced expression}. We use the following elementary property:
$\ell(ws_i)\in\{\ell(w)-1,\ell(w)+1\}$, $w\in\S$, $i=1,\dots,k-1$.
The elements of $\CC_{ki}=
\{p\in\S\colon p(m)<p(n)\text{ for } 1\le m<n\le i \text{ and }
i+1\le m<n\le k\}$ are called {\em shuffle}\/ permutations. Each $p\in\S$
admits a unique representation $p=p_1p_2p_3$ where $p_1\in\CC_{ki}$,
and $p_2\in \S$ (resp.\  $p_3\in \S$) leaves  $i+1,\dots,k$ 
(resp.\  $1,\dots ,i$) fixed. Moreover $\ell(p)=
\ell(p_1)+\ell(p_2)+\ell(p_3)$.
\\
{\em Lift into braids.} 
\/Artin's braid group $\B_k$ has generators $b_1,\dots,b_{k-1}$
and  defining relations 
\begin{xalignat}{2}
b_ib_{i+1}b_{i}&=b_{i+1}b_{i}b_{i+1},&\,&i=1,\dots,k-2,     \label{zopf1}\\
b_ib_j&=b_jb_i,      &\,&|i-j|\ge 2.                        \label{zopf2}
\end{xalignat}
The map $b_i\mapsto s_i$ defines  a natural
projection of $\B_k$ onto $\S$. For a reduced expression
$w=s_{i_1}\cdots s_{i_r}$  define  $b_w=b_{i_1}\cdots b_{i_r}$.
By \rf[zopf1] and \rf[zopf2], the
definition of $b_w$ does not depend  on the choice of the reduced 
expression $s_{i_1}\cdots s_{i_r}$. Obviously,
$
b_{vw}=b_vb_w$ for $\ell(vw)=\ell(v)+\ell(w).
$
This equation in particular applies to the shuffle decomposition $w=p_1p_2p_3$
of $w$:
$b_w=b_{p_1}b_{p_2}b_{p_3}$. Define the antisymmetrizer and
shuffle sums in the group algebra $\C\B_k$  as follows:
\begin{equation}
\label{symdefi}
\AAAA_k=\ts\sum_{w\in \S}(-1)^{\ell(w)} b_w, \qquad
\AAAA_{ki}=\ts\sum_{w\in \CC_{ki}}(-1)^{\ell(w)}b_w.
\end{equation}
By the shuffle decomposition we obtain for $1\le i<k$, $\AAAA_1=1$:
\begin{equation}
\AAAA_k=\AAAA_{ki}(\AAAA_i\ot \AAAA_{k-i}).    \label{sym}
\end{equation}
Let $\wo$ denote the longest word in $\S$. This permutation
maps $(1,\dots,k)$ into $(k,\dots,1)$.  
For a reduced expression $w=s_{n_1}\cdots s_{n_r}$ define 
$\wbar=s_{\ov{n_1}}\cdots s_{\ov{n_r}}$, where $\ov{n}=k-n$. Obviously,
$\wbar$ does not depend on the choice of the reduced expression. Using 
$s_n\wo=\wo s_{k-n}$ and $\ell(\ov{w})=\ell(w)$ one checks
\begin{equation}\label{bwo} 
b_wb_\wo=b_\wo b_\wbar \quad \text{and}\quad\AAAA_k b_\wo=b_\wo\AAAA_k.
\end{equation}
Now we recall the construction of the
external algebra $\gdw$ due to Woronowicz \cite[Proposition\,3.1]{Wo2}. There exists  a unique isomorphism 
$\sig\colon\Gamm\ota\Gamm\to\Gamm\ota\Gamm$, of bicovariant bimodules
called the {\em braiding}\/   with $\sig(\alpha\ota\beta)=\beta\ota\alpha$,
$\alpha\in\Gamm_{\ll}$, $\beta\in\Gamm_{\rr}$. Moreover $\sig$ fulfils the
braid equation \rf[zopf1], i.\,e.\
$\sig_{1}\sig_{2}\sig_{1}=\sig_{2}\sig_{1}\sig_{2}$ in
$\Gamm\ota\Gamm\ota\Gamm$, where  $\sig_{1}=\sig\ot \id$ and  $\sig_{2}=\id
\ot\,\sig$. Consequently, the map $\gamma(b_i)=\sig_i$, where
$\sig_i=\id\ott\cdots\ott\sig\ott\cdots\id$ acts in position $(i,i+1)$ of
$\gten[k]$, can  be extended to an algebra homomorphism
$\sigg\colon\C\B_k\to\End_\AA(\Gamm^{\ot k})$. We briefly write
$\sig_w=\sigg(b_w)$, $w\in \S$. Let  $\ant=\sigg(\AAAA_k)$ and
$\ant[ki]=\sigg(\AAAA_{ki})$. Call $\ant$ {\em Woronowicz' antisymmetrizer}\/ of
$\gten[k]$. By \rf[sym], $\JJ=\bigoplus_{k\ge 2}\ker \ant$ is a two-sided ideal
and a bicovariant bimodule. Hence $\gdw=\gten/\JJ$ is an $\NO$-graded algebra
and a bicovariant bimodule over $\AA$. Since $\gten[k]$ is a free left
$\AA$-module we always identify the linear spaces $(\gten[k])_\ll$ and
$\Gamm_\ll\ott\cdots\ott\Gamm_\ll$ ($k$ factors). Since $\sig$ is a morphism
of left comodules, it maps $(\Gamm\ota\Gamm)_\ll$ into itself. The
corresponding matrix form is $\sig(\theta_i\ot\theta_j)=\sum_{m,n}\sig^{mn}_{ij}\theta_m\ot\theta_n$, $i,j\in K$, with
\begin{equation}\label{fv}
\sig_{ij}^{mn}=f^i_n(v^m_j).
\end{equation}
\\
Recall that an $\NO$-graded algebra $H=\bigoplus_{n\ge 0} H^n$ is called {\em
  $\NO$-graded super  Hopf algebra}\/ if the product in $H\ot H$ is given by
$(a\ot b)(c\ot d)=(-1)^{ij}ac\ot bd$, $b\in H^i,\,c\in H^j$, $a,\,d\in H,$ and
there are linear mappings $\Delta$, $\ve$, and $S$ of degree $0$ called
coproduct, counit and antipode, respectively,   such that the usual Hopf
algebra axioms are fulfilled. Let $\Gamm$ be  a bicovariant bimodule over
$\AA$. The Hopf algebra  structure of $\AA$ uniquely extends to an
$\NO$-graded super Hopf algebra structure on $\gten$ such that  for
$\omega\in\Gamm$: $\Delta(\omega)=\dl(\omega)+\dr(\omega)$, $\ve(\omega)=0$,
and  $S(\omega)=-\sum S(\omega_{(-1)})\omega_{(0)}S(\omega_{(1)})$, see
\cite[Proposition\,13.7]{KS}. The antipode is a graded anti-homomorphism i.\,e. 
\begin{equation}  \label{srr}
S(\rho_1\ota\rho_2)
=(-1)^{kn}S\rho_2\ota S\rho_1\quad\text{ for }\quad
\rho_1\in\gten[k],\, \rho_2\in\gten[n]. 
\end{equation}
Moreover  $\Js$  and $\JJ$ are Hopf ideals in  $\gten$, see  \cite[p.~489
before Proposition~13.9 and Proposition~13.10]{KS}, and so is $\Ju$, see
\cite[Theorem~14.8]{KS}.  Hence, $\gdu$, $\gds$, and $\gdw$ are  $\NO$-graded
super Hopf algebras. We write $\alpha\beta$ and $\gd[k]$ instead of
$\alpha\ota\beta$ and $\gten[k]$, respectively, when dealing with one of
these quotients  $\gten{/}J$. Let $B$ be an $\NO$-graded algebra,
$B=\bigoplus_{k\ge 0}B_k$, then the formal power series $P(B,t)=\sum_{k\ge
0}(\dim B_k)\,t^k$, is called {\em Poincar\'e series}\/ of $B$.
\\
\subsection*{Bicovariant Differential Calculus}
A {\em first order differential calculus}\/ over $\AA$ abbreviated  FODC  is an
$\AA$-bimodule $\Gamm$ with  a linear mapping  $\dd\colon\AA\to\Gamm$ that
satisfies the Leibniz rule $\dd(ab)=\dd a{\cdot} b+a{\cdot}\dd b$ for
$a,\,b\in\AA$,  and $\Gamm$ is the linear span of elements $a\dd b$ with
$a,b\in\AA$.
\\
A {\em differential graded algebra}\/ over $\AA$ is  an  $\NO$-graded algebra
$\gd=\bigoplus_{n\ge0}\Gamm^n$, $\Gamm^0=\AA$, with  a linear mapping
$\dd\colon\gd\to\gd$ of degree $1$ such that $\dd^2=0$, and $\dd$
satisfies the graded Leibniz rule $\dd(\rho_1\rho_2)=\dd\rho_1{\cdot}\rho_2
+(-1)^n \rho_1{\cdot}\dd \rho_2$, $\rho_1\in\Gamm^n$, $\rho_2\in\gd$. If in
addition $\Gamm^n=\AA{\cdot}\dd\AA\cdots\dd\AA$ ($n$ factors), $n\in\N$, call
$\gd$ a {\em differential calculus}. 
A differential calculus $\gd$ is called {\em bicovariant}\/ if there exist linear mappings
$\dl\colon\gd\to\AA\ot\gd$ and
$\dr\colon\gd\to\gd\ot\AA$ of degree $0$ with
$\dl{\uhr}\AA=\dr{\uhr}\AA=\Delta$ such that
\\
(i) $(\gd,\dl,\dr)$ is a bicovariant
bimodule, and
\\
(ii) $\dl(\dd \rho)=(\id\ot\dd)\dl(\rho)$, and $\dr(\dd
\rho)=(\dd\ot\id)\dr(\rho)$ for $\rho\in\gd$. 
\\
A differential calculus is called {\em inner}\/ if there exists an element  $\rho\in\Gamm^1$ such that
$\dd\rho_n=\rho\rho_n-(-1)^{n}\rho_n\rho$, $\rho_n\in\Gamm^n$. 
\\
A FODC $\Gamm$ is called {\em bicovariant}\/ (resp.\  {\em inner}\/) if the
differential calculus $\AA\oplus\Gamm$ is bicovariant (resp.\  inner).
Let $(\Gamm_i,\dd_i)$, $i=1,2$, be differential graded algebras over
$\AA$. Then $(\Gamm_1\ot\Gamm_2,\dd_\ot)$ becomes a differential graded
algebra  over $\AA$ if the product  in $\Gamm_1\ot\Gamm_2$ is defined  by
$(\omega_1\ot\omega_2)(\rho_1\ot\rho_2)=(-1)^{ij}\omega_1\rho_1\ot\omega_2\rho_2$,
$\omega_2\in\Gamm_2^{i}$, $\rho_1\in\Gamm_1^{j}$, $\omega_1\in\Gamm_1$,
$\rho_2\in\Gamm_2$, and the differential $\dd_\ot$ is given by 
\begin{equation}
  \nn
\dd_\ot(\omega_1\ot\omega_2)=\dd_1\omega_1\ot\omega_2+
(-1)^{i}\omega_1\ot\dd_2\omega_2,
\quad\omega_1\in\Gamm_1^{i},\,\omega_2\in\Gamm_2.
\end{equation}
A  {\em  differential Hopf algebra}\/ is both a differential graded algebra with
differentiation $\dd$ and an $\NO$-graded super Hopf algebra  with coproduct $\Delta$ satisfying the condition $\dd_\ot\Delta=\Delta\dd$. 
\\
Our main objects $\gdw$, $\gds$, and $\gdu$ are  differential  Hopf algebras, see \cite[Theorem\,14.17 and Theorem\,14.18]{KS} for details.
\\
{\em Associated right ideal, the mappings $\om$ and $\SS$.}\/  According to
\cite[Theorems\,1.5 and 1.8]{Wo2}  there is a one-to-one correspondence
between bicovariant first order differential calculi $\Gamm$ over $\AA$ and
$\adr$-invariant right ideals $\RR$  in $\ker\ve$ given by
$\RR=\ker\ve\cap\ker\om$. A crucial role play the two mappings
$\omega\colon\AA\to\gl$, $\om(a)=\sum Sa_{(1)}\dd a_{(2)}$, and
$\SS\colon\AA\to\gl\ot\gl$, $\SS(a)= \sum \om(a_{(1)})\ot\om(a_{(2)})$. Note
that $\Ju$ is the ideal in $\gten$ generated by  $\SS(\RR)$. Both $\om$ and
$\SS$ intertwine the right adjoint coaction $\adr$ with the right coaction
$\dr$ on $\Gamm_\ll$ and $\Gamm_\ll\ot\Gamm_\ll$, resp. The mappings $\om$ and
$\SS$ are coupled by the Maurer-Cartan equation, see
\cite[Proposition\,14.13]{KS}:
\begin{equation}
\label{MC}
\dd\om(a)=-\ts\sum_{(a)}\om(a_{(1)})\om(a_{(2)}),\,a\in\AA.
\end{equation}
\subsection*{Bicovariant Differential Calculus on quantized simple Lie groups}
Throughout  the deformation parameter $q$ is assumed to be a {\em
transcendental}\/ number, and $N\ge2$. Let $\AA$ be  the
bialgebra $\AA(R)$ or one of the Hopf algebras  $\oglqn$ and $\oslqn$ as
defined in \cite[Subsect.\,1.3]{FRT}.  Recall that $R$ denotes the
 complex invertible $N^2\times N^2$-matrix $R^{ab}_{rs}=\Rda[rs]{ba}$, where
$\Rda{}$ is given in  \rf[rda]. The  $N^2$ generators of $\AA$ are denoted
by  $u^i_j$, $i,\,j=1,\dots,N$, and we call $u=(u^i_j)_{i,j=1,\dots,N}$ the
{\em fundamental matrix corepresentation}. As an algebra $\oglqn$ is
generated by the elements  $u^i_j$ and an additional central element
$\dett$ satisfying $\DD\dett=\dett\DD=1$, where
$\DD=\sum_{p\in\S[N]}(-q)^{\ell(p)}u^1_{p(1)}\cdots u^N_{p(N)}$ denotes the
quantum determinant. The element  $U=\sum_{i=1}^Nq^{-2i}u^i_i$  is called
quantum trace.  We abbreviate $\QM=q-q^{-1}$ and $\QP=q+q^{-1}$. The $\Rda{}$-matrix is given by
\begin{equation}
\label{rda}
\Rda[rs]{ab}=q^{\delta_{ab}}\delta_{as}\delta_{br}
+\QM h(b-a)\delta_{ar}\delta_{bs},
\end{equation}
$a,\,b,\,r,\,s=1,\dots,N$, where $h$ denotes the Heaviside symbol  $h(x)=1$
for $x>0$ and $h(x)=0$ for $x\le0$. The matrix $\Rda{}$ can be written as
$\Rda{}=qP_+-\qm P_-$, where $P_\pm=\QP^{-1}(q^{\mp 1}I\pm\Rda{})$ are
projections. 
\\
We follow the method of \cite{J} and \cite{CSWW} to construct bicovariant FODC
on quantizations of simple Lie groups. For  a nonzero complex number
$x\in\C^\times$ let  $\ell^\pm_x=((\ell_x^\pm)^i_j)$ be the $N\times N$-matrix of linear functionals $(\ell_x^\pm)^i_j$ on $\AA=\oglqn$ as defined in
in \cite[Sect.\,2]{FRT}. Recall that $\ell^\pm_x$ is uniquely determined by
${\ell^{\pm}_{x}}^i_j(u^m_n)=x^{\mp 1}(\Rda{\pm 1})_{nj}^{im}$ and the
property that  $\ell^\pm_x\colon\AA\to M_N(\C)$ is a unital algebra
homomorphism. For $\AA=\oslqn$ we must assume $x^N=q$. Define the bicovariant
bimodules $\gtz$,  $\tau\in\{+,-\}$, where $z\in\C^\times$ is arbitrary  in
case $\AA=\oglqn$ and $z^N=q^{-2\tau}$ in case $\AA=\oslqn$:
\begin{equation}          \label{gpmz}
\Gamm_{+,z}=(u^\fettc\ot u,\ell^+_x\ot\ell^{-,\fettc}_y),   \quad xy=z^{-1},\quad
\Gamm_{-,z}=(u^\fettc\ot u, \ell^-_x\ot \ell^{+,\fettc}_y),\quad xy=z.
\end{equation}
The structure of $\gpmz$ can easily be described as follows. There exists a
basis $\{\theta^i_j:i,\,j=1,\dots,N\}$ of $(\gpmz)_\ell$ such that  the right
action and the right coaction  are given by
\begin{align}
\label{oma}
\theta^i_j a&=\sum_{m,n}({\ell^\pm}^i_mS{\ell^\mp}^n_j\ast a)\theta^m_n,\quad a\in \AA  \\
\intertext{and}           
\dr \theta^i_j&=\sum_{m,n}\theta^m_n\ot(u^\fettc)^m_iu^n_j,\quad i,\,j=1,\dots,N.\label{right}
\end{align}
The element $\theta=\sum_iq^{-2i}\theta^i_i$ is the unique up to scalars
bi-invariant element. Defining 
\begin{equation}
\label{inner}
\dd a=\theta a-a\theta
\end{equation}
for $a\in \AA$, $(\gpmz,\dd)$ becomes a bicovariant FODC over $\AA$. The
braiding $\sig$ of $\gpmz$ can be obtained as follows. Inserting
$v=u^\fettc\ott u$ and $f=\ell_x^{\pm}\ot\ell_y^{\mp,\fettc}$ into \rf[fv] yields
\begin{equation}\label{sig}
\sig_\pm=\Rre[23]^{\pm}\Rch[12]{\mp1}
\Rda[34]{\pm1}(\Rre[23]^{\pm})^{-1},
\end{equation}
where the complex $N^2{\times}N^2$-matrices $\Rre^\pm$ and $\Rch{}$ are
defined as follows: $\Rch[rs]{ab}=\Rda[ba]{sr}$ and
$(\Rre^\pm)_{rs}^{ab}=(\Rda{\pm1})^{ra}_{sb}$.  The indices $12$, $23$, and
$34$ indicate that the corresponding matrices act at positions  $1,2$; $2,3$; and $3,4$ of the tensor product $V^{\ot 4}$, $V=\C^N$. 
Note that the braiding $\sig$ of $\gtz$ does not  depend on $z$. Let
$r_+=\qm\QM$ and $r_-=-q^{-2N-1}\QM$. We often write  $\sig_\tau$, $r_\tau$,
$\Rda{\tau}$, and $\Rre^{\tau}$ instead of $\sig_\pm$, $r_\pm$, $\Rda{\pm 1}$,
and $\Rre^{\pm}$, respectively. Set $\fr{s}=q^{-2}+q^{-4}+\cdots +q^{-2N}$,
$\fr[\tau]{r}=\fr{s}+r_\tau$, and $\ntz=z\fr[\tau]{r}-\fr{s}$. Corresponding
to $\gtz$ the parameter value $z$ is called {\em regular}\/ if $\ntz\ne0$. For
regular $z$  the differentials $\dd u^i_j$ generate $\gtz$ as a  left
$\AA$-module. Equivalently, the linear space
$\langle\om_{ij}:i,\,j=1,\dots,N\rangle$ of  Maurer-Cartan forms
$\om_{ij}=\om(u^i_j)$ is $N^2$-dimensional. The  right ideal $\RR_\tz$ then satisfies
\begin{equation}\label{lin}
\RR_\tz\oplus\langle u^i_j:i,j=1,\dots,N\rangle\oplus\C=\AA,
\end{equation}
see \cite[Lemma\,1.5]{SS2} and \cite[Theorem\,2.1]{SS3}. For $\AA=\oslqn$,
$\Gamm=\gtz$, $z^N=q^{-2\tau}$, the  condition $\ntz\ne0$ is trivially
fulfilled since $q$ is  transcendental. Transformation formulas between
$\theta^i_j$, $\theta$, and $\om_{ij}=\om(u^i_j)$, $\vt
=\sum_iq^{-2i}\om^i_i$, are given as follows. Using matrix notation and leg
numbering for the matrices $\vO=(\om_{ij})$, $\vN=(\theta^i_j)$, and $u$,
where $\vN_{1}=\vN\ot I$, $\vN_{2}=I\ot \vN$ and $u_{2}=I\ot u$, we have
\begin{equation}
\label{e-trans}
\vO    =zr_\tau\vN +(z-1)\theta I\quad
\text{and}\quad\vt=\ntz\theta.
\end{equation} 
Inserting  $a=(u^k_l)$ into \rf[oma] gives
\begin{xalignat}{2}     
\theta^i_j\rac u^k_l&=z\sum_{x,m,n}
(\Rda{\tau})^{ik}_{xm}\theta^m_n(\Rda{\tau})^{xn}_{jl},& \theta\rac
u^i_j&=zr_\tau\theta^i_j+z\theta\delta_{ij},\nn
\\
\intertext{or in matrix notation}
\vN_{1}\rac u_{2}&=z \Rda{\tau}\vN_{2}\Rda{\tau},&
\theta\rac u  &=z r_\tau\vN +z\theta I. \label{rac}
\end{xalignat}
\section{Main results}\label{result}
\begin{thm}     \label{left}
Let $\AA$ be one of the Hopf algebras $\oslqn$ or $\oglqn$ and let $\Gamm$ be
one of the $N^2$-dimensional bicovariant first order differential calculi
$\gtz$, $\tau\in\{+,-\}$, $z\in\C^\times$,  over $\AA$. Suppose $q$ to be
transcendental.  Let $\gdw$ denote Woronowicz' external algebra and let $\gdwl$ denote its subalgebra of left-invariant forms.
\\
Then the Poincar\' e series of $\gdwl$ is
\begin{equation*}
P(\gdwl,t)=(1+t)^{N^2},
\end{equation*}
i.\,e.  $\dim (\gdwl[k])=\binom{N^2}{k},\,k\in\NO$. In particular, there
exists a unique up to scalars left invariant form of maximal degree $N^2$.
\end{thm}
\begin{thm}\label{bi}
Let $\AA$,  $\Gamm$, and $q$ be as in Theorem\,\ref{left}. Let $\gdwi$ denote
the subalgebra of $\gdw$ consisting of all bi-invariant forms. Then
\smallskip\\ 
{\rm (i)} $P(\gdwi,t)=(1+t)(1+t^3)\cdots(1+t^{2N-1})$.
\\
{\rm (ii)} $\omega_k\land\omega_n=(-1)^{kn}\,\omega_n\land\omega_k,$ for $\omega_k\in {\gdwi[k]}$ and $\omega_n\in{\gdwi[n]}$.
\\
{\rm (iii)} $\dd \omega=0$ for $\om\in\gdwi$. 
\\ 
{\rm (iv)}  Different bi-invariant forms represent different de Rham
cohomology classes. 
\end{thm}
\begin{remark} By (i) the dimension $\dim(\gdwi[k])$ is equal to the number of
partitions of $k$ into a sum of pairwise different positive odd integers less
than $2N$. It is also equal to the number of  symmetric partitions $\lam$ of $k$, $\lam=\lam'$, with  $\lambda_1\le N$. The left-invariant form of maximal degree $N^2$ is bi-invariant.
\end{remark}
\begin{thm}  \label{usw}
Let $\AA$, $\Gamm$, and  $q$  be as in Theorem\,\ref{left}. Let $\gdu$ and
$\gds$ denote the universal differential calculus and the second
antisymmetrizer  differential calculus over $\Gamm$, respectively. Then we
have the following isomorphisms of differential Hopf algebras. 
\smallskip\\ 
{\rm (i)} $\gds  \cong\gdw$.  
\\ 
{\rm (ii)} $\gdu\cong\gds$  for $N\ge3$, $\ntz\ne0$  and for $N=2$, $\tau=+$,
and
\\
\quad$z\ne\qm,-\qm,(q^2+1)(q^4+1)^{-1}$.
\\
{\rm (iii)}  $\gds\cong\gdu/(\theta^2)$ for  $N=2$, $\tau=+$, and $z^2=\qm[2]$. Moreover, the universal differential calculus $\gdu$ is inner and $\theta^2$ is central.
\end{thm}
\begin{remark} For $\AA=\mathcal{O}(GL_q(n|m))$ and $\Gamm=\Gamm_{\pm,1}$, (ii)
was proved in \cite{LS}. \\
Since the differential calculi  $\Gamm_{+,z}$  and  $\Gamm_{-,z}$
are isomorphic for $N=2$ it suffices to consider $\Gamm_{+,z}$. In case of the quantum
group $\slqn[2]$, the parameter $z=q^{-1}$ (resp.\  $z=-q^{-1}$) corresponds
to the $4D_+$-calculus  (resp.\  $4D_-$-calculus),  defined in  \cite{Wo2}. \\
The calculus $\gds$  is also inner since $\theta^2=0$, see Lemma\,\ref{l-inner}\,(ii) below.
\end{remark}
\section{Iwahori-Hecke algebra and Antisymmetrizer}\label{s-hecke}
In this subsection we shall give the definition of the Iwahori-Hecke algebra
$\h$ over $\C$ and we shall list some of their important properties and facts
including an explicit formula for the central idempotents. Since $q$ is not a
root of unity everything goes through in exactly the same fashion as in the
case of the base field $\C(q)$. We take the definition from
\cite[Sects.~7.1 and 7.4]{Hum}.
\\
There is a unique structure of an associative unital algebra on  the  vector
space with basis $\{T_w\colon w\in\S\}$ and $T_1=1$ such that for all
$s\in\{s_1,\dots,s_{k-1}\}$ and $w\in\S$ 
\begin{alignat}{2}
T_sT_w&=T_{sw}&\text{ if }\quad\ell(sw)>\ell(w),              \label{h1}
\\
T_sT_w&=(q-q^{-1})T_w+T_{sw}&\text{ if }\quad\ell(sw)<\ell(w).    \label{h2}
\end{alignat}
This algebra is called {\em Iwahori-Hecke algebra of type}\/ ${\rm A}_{k-1}$
and is denoted by $\h$. 
Instead of \rf[h2] one often takes $T_sT_w=(q-1)T_{w}+ qT_{sw}$. The present
form is more useful when dealing with quantum groups. We briefly write $T_i$
for $T_{s_i}$. The relations \rf[zopf1], \rf[zopf2] with $T_i$ instead of
$b_i$ and 
\begin{equation}  \label{T2}
T_i^2=\QM T_i +1,\quad i=1,\dots,k-1,    
\end{equation}
are equivalent  to \rf[h1] and \rf[h2].
\\
We shall adopt the notations in \cite{Mac}  for partitions, compositions,
characters, and idempotents.  We call $\lam=(\lam_1,\lam_2,\dots )$  a {\it
  partition}\/ of $k$ if  $\lam_1+\lam_2+ \dots=k$ and
$\lam_1\ge\lam_2\ge\cdots\ge0$ are integers. We briefly write $\lam\vdash
k$. Sometimes we use the notation $\lam=(k^{m_k}(k-1)^{m_{k-1}}\cdots
1^{m_1})$ where $m_i$ denotes the number of parts of $\lam$ equal to $i$. The
{\it conjugate}\/ of a partition $\lam$ is a partition $\lam^\prime$ whose
diagram is the transpose of the diagram of $\lam$. Hence $\lam_i^\prime$ is the
number of boxes in the $i$th column of $\lam$, or equivalently
$\lam_i^\prime={\rm Card}\{j\colon \lam_j\ge i\}$. 
\\
Irreducible representations of $\h$  are labelled by partitions $\lam$ of
$k$. The  dimension of the irreducible representation labelled by $\lam$ is
$d_\lam=\frac{k!}{h(\lam)}$. Here $h(\lam)=\prod_{x\in\lam}h(x)$ is
the product of {\em hook-lengths}\/ $h(x)=\lam_i+\lam'_j-i-j+1$, where
$x=(i,j)$ is a box of $\lam$ in position $(i,j)$, see \cite[I.\,6,
Example\,2\,(a) and I.\,7 (7.6)]{Mac}. The corresponding irreducible character
is denoted by $\chi^\lam$.  A {\it partition of unity}\/ is a set of minimal
(primitive) idempotents $\{p^i_\lam\}$, $\lam\vdash k$, $1\le i\le d_\lam$,
such that $p^i_\lam p^j_\mu=\delta_{\lam\mu}\delta_{ij}$. The element
$z_\lam=\sum_{i=1}^{d_\lam}p_\lam^i$ is the  minimal {\em central}\/ idempotent
associated to  $\lam$. We have $1=\sum_\lam z_\lam$. For the remainder of this
section fix positive integers $k$ and $N$. An {\it $N$-composition}\/ of $k$ is
a sequence of nonnegative integers $c=(c_1,\dots,c_N)$, $c_i\ge 0$, such that
$c_1+c_2+\dots +c_N=k$. It is denoted by $c\models k$. For any $N$-composition
$c=(c_1,\dots,c_N)$ of $k$ define $x^c=x_1^{c_1}\cdots x_N^{c_N}$, where
$x_1,\dots,x_N$ are commuting variables. Let $\{e_1,\dots,e_N\}$ be a basis of
the vector space $V=\C^N$. For $v=e_{i_1}\ott e_{i_2}\ott\cdots\ott e_{i_k}\in
V^{\ot k}$  define the {\it content} $c(v)$, $c(v)=(c_1,\dots,c_N)$, where
$c_j$ is the number of factors $e_j$ in $v$. Obviously, $c(v)\models k$. For
$c\models k$ define the projection operator $E_c$ of $V^{\ot k}$ by
$E_c(v)=\delta_{c,c(v)}v$ for $v=e_{i_1}\ott e_{i_2}\ott\cdots\ott
e_{i_k}$. One easily verifies $E_cE_{c'}=\delta_{c,c'}E_c$, $c,c'\models k$,
and $I=\sum_{c\models k}E_c$.
\\
We introduce two important involutions $^\ast$ and $^\prime$ of $\h$. They are
defined on generators by $T_s^\ast=T_s$ and $T_s^\prime=-T_s^{-1}$,
respectively. Both mappings can be extended to antihomomorphisms of $\h$. Note
that $T_w^\ast=T_{w^{-1}}$ and $T_w^\prime=(-1)^{\ell(w)}T_w^{-1}$.
\begin{lemma}\label {invzentral}
Let $\lam\vdash k$ be a partition of $k$. Then $z_\lam^\prime=z_{\lam^\prime}$ and $z_\lam^\ast=z_\lam$.
\end{lemma}
\begin{proof}
(a) The central idempotent $z_{(k)}$ (resp.\  $z_{(1^k)}$) is completely
characterized by the property $T_iz_{(k)}=qz_{(k)}$ (resp.\
$T_iz_{(1^k)}=-q^{-1}z_{(1^k)}$), $i=1,\dots,k-1$. Since any involution of a
central idempotent is again a central idempotent,
$T_iz_{(k)}^\prime=-(z_{(k)}T_i^{-1})^\prime=-q^{-1}z_{(k)}^\prime$,
$i=1,\dots,k-1$, implies $z_{(k)}^\prime=z_{(1^k)}$. Since the
anti-automorphism $^\prime$ is of order two, $z_{(1^k)}^\prime=z_{(k)}$. The proof of $z_{(k)}^\ast=z_{(k)}$ and $z_{(1^k)}^\ast=z_{(1^k)}$ is similar.
\\ 
(b) Suppose we are given an admissible Young tableau with diagram
$\lam=(\lam_1,\dots,\lam_r)$. Let
$\lam^\prime=(\lam_1^\prime,\dots,\lam_s^\prime)$ be the transposed diagram.
In \cite{Gy} there is constructed a minimal idempotent
$p_\lam=h_-e_-h_-^{-1}h_+e_+h_+^{-1}$, where
$e_-=A_{\lam_1^\prime}\times\cdots \times A_{\lam_s^\prime}$,
$e_+=S_{\lam_1}\times\cdots \times S_{\lam_r}$. Here $A_j=z_{(1^j)}$ and
$S_j=z_{(j)}$ denote the $j$th antisymmetrizer and the $j$th symmetrizer in
$\h[j]$, respectively, and we use the natural embeddings
$\h[n_1]\times\cdots\times\h[n_r]\to\h$ for  $n\models k$. The element $h_\pm$
depends on the entries  of  the Young tableau. By (a),
$e_-^\prime=S_{\lam_1^\prime}\times\cdots \times S_{\lam_s^\prime}$ and
$e_+^\prime=A_{\lam_1}\times\cdots \times A_{\lam_r}$. Since $p_\lam$ is a
minimal subidempotent of $z_\lam$, $p_\lam^\prime$ is a minimal subidempotent
of $z_{\lam'}$.  Since $z_\lam=\sum_{i=1}^{d_\lam}p_\lam^i$ and
$(p_\lam^i)^\prime(p_\lam^j)^\prime=\delta_{ij}(p_\lam^i)^\prime$,
$z_\lam^\prime$ is a subidempotent of $z_{\lam'}$. Moreover, $z_\lam^\prime$
is non-zero and a  minimal central idempotent. Hence,
$z_\lam^\prime=z_{\lam'}$. The proof for  $^\ast$ is similar. 
\end{proof}
Next define  representations $\vr$ and $\vr_\fettc$ of $\h$ on $V^{\ot
k}$. The  action of the generator $T_n$ on a basis element $v=e_{i_1}\ot
e_{i_2}\ot\cdots \ot e_{i_k}$  is given as follows. Let $s_n$ denote the
flip operator in $V^{\ot k}$ that interchanges the $n$th and $(n+1)$st component of the tensor product. Let
\begin{equation*}
\vr(T_n) v=\begin{cases}
                qv          &\text{ if }\quad i_n=i_{n+1},\\
                \QM v+s_n(v)&\text{ if }\quad i_n<i_{n+1},\\
               s_n(v)        &\text{ if }\quad i_n>i_{n+1}.
         \end{cases}
\end{equation*}
The  representation $\vr_\fettc$ is given by
$\vr_\fettc(T_n)=s_n\vr(T_n)s_n$. Using matrix notation and leg numbering,
$\vr(T_n)=\Rda[n, n+1]{}$. The representation $\vr_\fettc$ associates to  a
single generator $T$ the matrix $\check R=(\check R^{ij}_{kl})=(\hat
R^{lk}_{ji})$. Both representations $\vr$ and $\vr_\fettc$ of $\h$ commute
with $E_c$ since they do not change the content of the  tensor
$e_{i_1}\ot\cdots\ot e_{i_k}$.
\\
We  introduce a trace functional $\TR$ on $\h$ that takes values in the
algebra of polynomials $\C[x_1,\dots,x_N]$. The following Proposition is due
to Ram, see \cite[Lemma\,3.5, Lemma\,3.7 and  Theorem 3.8]{R}.
\begin{prp}\label{p-heckespur} 
{\rm (i)} The mapping $\TR\colon \h\to\C[x_1,\dots,x_N]$
\begin{equation}\label{d-tau}
\TR(h)=\sum_{c\models k}x^c\tr(E_c\vr(h))
\end{equation}
defines a trace functional, i.\,e.\  $\TR(hg)=\TR(gh)$ for all $h,g\in\h$.
\\ {\rm (ii)}\quad Let  $p\in \h$ be idempotent. Then $\TR(p)$ does not depend
on $q$. In particular 
\begin{equation}
\TR(p_\lam^i)=s_\lam(x_1,\dots,x_N), \label{spurpl}
\end{equation}
where $\lam\vdash k$ is a partition of $k$, $\{p_\lam^i\}$ is a partition of
unity, and $s_\lam$ denotes the Schur function, cf. {\rm\,\cite[I\,(3.1), p.~40]{Mac}}.
\\
{\rm (iii)} For  $h\in\h$ 
\begin{equation}\label{spurchar} 
\begin{split}
\TR(h)&=\sum_{\lam\vdash k}\chi^\lam(h)s_\lam(x_1,\dots,x_N),\\
\TR(hz_\lam)&=\chi^\lam(h)s_\lam(x_1,\dots,x_N).
\end{split}
\end{equation}
\end{prp}
\begin{remark}\label{hecke}
The proof of \cite[Lemma 3.5]{R} shows that replacing $\vr$ by $\vr_\fettc$
in \rf[d-tau] does not change $\TR$. 
\end{remark}
For the group algebra $\C\S$ we have $z_\lam=\frac{d_\lam}{k!}\sum_{w\in
\S}\chi^\lam(w)w^{-1}$. We will show now that a similar formula holds for
$\h$. The existence of an associative, symmetric, and non-degenerate bilinear
form on $\h$ is essential for the proof of this formula. Let $\{T_w: w\in\S\}$
and $\{f^w:w\in\S\}$ be dual bases of $\h$  and its dual vector space, respectively. Let $f^0$ denote the coordinate functional corresponding to the unit element $1\in\h$.
\begin{lemma}\label{l-bilinear} 
Let $g,h\in\h$. Then $\langle g,h\rangle=f^0(gh)$ defines a symmetric,
associative, and non-degenerate bilinear form on $\h$. Moreover
$\langle T_v,T_w\rangle= 
\begin{cases}1,&\text{ if } vw=1\\
             0 &\text{ otherwise}.
\end{cases} 
$\hfill $(\ast)$
\end{lemma}
\begin{proof}
For the right hand side of $(\ast)$ we use the Kronecker  symbol $\delta_{vw,1}$. The
associativity of the pairing follows from the definition. Suppose for a moment
that formula $(\ast)$ is already proved. Then $\pair$ is symmetric since
$vw=1$ if and only if $wv=1$, and $f^0$ is linear. The pairing is
non-degenerate since $\{T_w:w\in\S\}$ and $\{T_{w^{-1}}:w\in\S\}$ are
orthogonal with respect to $\pair$ bases. We prove $(\ast)$ by induction on
the length of $w$. By definition $\langle T_v,1\rangle=f^0(T_v)=\delta_{v,1}$,
and $(\ast)$ holds for $w=1$. Suppose  it  is true for all $v,w'\in\S$ with
$\ell(w')\le\ell(w)$. Let $w=sw'$, $\ell(w)=\ell(w')+1$, and $v$ be arbitrary. 
\\
{\em Case 1}. $\ell(vs)=\ell(v)+1$. By \rf[h1] and induction assumption $\langle
T_v,T_w\rangle=f^0(T_vT_sT_{w'})=f^0(T_{vs}T_{w'})=\langle
T_{vs},T_{w'}\rangle=\delta_{vsw',1}=\delta_{vw,1}$. 
\\
{\em Case 2}. $\ell(vs)=\ell(v)-1$. By \rf[h2] and induction assumption we
have $\langle T_v,T_w\rangle=f^0(T_vT_sT_{w'})=f^0\bigl((\QM
T_{v}+T_{vs})T_{w'}\bigr)=\QM\langle T_{v},T_{w'}\rangle+\langle
T_{vs},T_{w'}\rangle=\QM\langle T_v,T_{w'}\rangle+\delta_{vw,1}$. We have to
prove that $\langle T_v,T_{w'}\rangle=0$ or equivalently, by induction
assumption, $vw'\ne 1$. Assume to the contrary $vw'=1$. Then
$\ell(v)=\ell(w')$. Since also $vssw'=1$, $\ell(vs)=\ell(sw')=\ell(w')+1=\ell(v)+1$; that contradicts our assumption of case 2. Hence $vw'\ne 1$ and the proof is complete.
\end{proof}
Since there exists an associative, symmetric, and non-degenerate bilinear form
on $\h$, Proposition~(9.17), p.\,204 in  \cite{CR} applies to our
situation. Namely, the central idempotent $z_\lam$ is given by 
\begin{equation}                                           \label{e-zlam}
z_\lam={\ts\frac{d_\lam}{t_\lam}} \sum_{w\in\S}\chi^\lam(T_w)T_{w^{-1}},
\end{equation} 
where $t_\lam=\sum_{w\in\S}\chi^\lam(T_w)\chi^\lam(T_{w^{-1}})$ is nonzero. Moreover
\begin{equation}\label{e-tlam}
t_\lam=d_\lam^{-1}\chi^\lam(t)\quad\text{ with  }\quad 
t=\sum_{w\in\S}T_wT_{w^{-1}}.
\end{equation}
\begin{lemma}\label{l-tcentral}
The element $t=\sum_{w\in\S}T_wT_{w^{-1}}$ is central in $\h$. We have
$t=\sum_{\lam\vdash k} t_\lam z_\lam$. 
\end{lemma}
\begin{proof}
Fix $s=s_i$ for some $i\in\{1,\dots,k-1\}$ and set
$L=\{w\in\S:\ell(sw)>\ell(w)\}$. Then $\S$ is the disjoint union of $L$ and
$sL$. Abbreviating $\til=T_wT_{w^{-1}}$, by \rf[h1] we have
$\wt{T}_{sw}=T_s\wt{T}_wT_s$, $w\in L$. Moreover, by \rf[T2] 
\begin{equation}\nn
\begin{split} 
T_stT_s&=\sum_{w\in L\cup sL}T_s\til T_s=\sum_{w\in L}
(T_s\til T_s+ T_s^2\til T_s^2)
=\sum_{w\in L} (T_s\til T_s+T_s^2\til(\QM T_s+1))
\\
 &=\sum_{w\in L} (T_s\til T_s+(\QM T_s+1)\til +\QM T_sT_s\til T_s)
 \\
 &= \sum_{w\in L}(T_s\til T_s+\til)+
 \QM T_s\sum_{w\in L} (\til+T_s\til T_s)=t+\QM T_s t. 
\end{split}
\end{equation} 
Multiplying the preceding equation  from the left by $T_s^{-1}$ and using
$T_s^{-1} =T_s-\QM$, we get $tT_s=T_s t$. Hence, $t$ is central. Applying $\chi^\lam$ to the ansatz $t=\sum_{\mu\vdash k}\alpha_\mu z_\mu$,
using $\chi^\lam(z_\mu)=\delta_{\lam\mu}d_\lam$, and \rf[e-tlam], we
obtain $\alpha_\lam=t_\lam$.
\end{proof}
In our studies the algebra $\H=\h\ot\h$ plays the important role. The main
object is the antisymmetrizer $\AAA_k$ in $\H$. We define it as follows. Let
the map $\bsig$ be defined on the generators of the braid group  by
$\bsig(b_i)=T_i^{-1}\ot T_i$. Since  the elements of $\{\bsig(b_i)\}$ satisfy
the braid equation \rf[zopf1] and \rf[zopf2], $\bsig$ uniquely extends to a homomorphism
of the group algebra $\C\B_k$ to $\H$. One easily checks that
$\bsig_w:=\bsig(b_w)=(-1)^{\ell(w)}T_w^{\ast\prime}\ot T_w$, $w\in \S$. Let
$\AAA_k$ denote the image of the braid group antisymmetrizer \rf[symdefi]: 
\begin{equation}                              \label{d-anti}
\AAA_k=\bsig(\AAAA_k)=\sum_{w\in\S}T_w^\prime\ot T_w^\ast.
\end{equation}
Similarly, define $\AAA_{ki}=\bsig(\AAAA_{ki})$ for all $i<k$. The basic
result of this section is 
\begin{prp}\label{spectrum} The
spectral decomposition of the antisymmetrizer $\AAA_k\in \H$ is 
\begin{align}\AAA_k&=\sum_{\lam\vdash k}t_\lam \pi_\lam, \label{spec}
\\
\intertext{where}
\pi_\lam&={t_\lam}^{-1}\AAA_k(z_\lam^\prime\ot z_\lam)\label{pia}
\end{align} are mutually orthogonal idempotents in $\H$. Moreover
\begin{equation}\label{ttPi}
(\TR\ot\TR)\pi_\lam=s_{\lam'}(x_1,\dots,x_N)s_\lam(x_1,\dots,x_N).
\end{equation}
\end{prp}
\begin{proof}
The following formula is essential for the proof of $\pi_\lam^2=\pi_\lam$: 
\begin{equation}\label{relianti}
\AAA_k(1\ot h)=\AAA_k(h^\prime\ot 1)
\end{equation}
for $h\in\h$. We prove \rf[relianti] in three steps. (a) $k=2$ and
$h=T=T_1$. By \rf[T2] and $T^{-2}=1-\QM T^{-1}$ we have
\begin{align*}
\AAA_2(1\ot T)&=(1\ot 1-T^{-1}\ot T)(1\ot T)=1\ot T-T^{-1}\ot(\QM T+1)
\\
&=
(1-\QM T^{-1})\ot T-T^{-1}\ot 1=(-T^{-1}\ot T+1\ot 1)(-T^{-1}\ot 1)=
\AAA_2(T'\ot 1).
\end{align*}
(b) Now let $k>2$. By \rf[sym],
$\AAA_k=\AAA_{k,i+1}(\AAA_{i+1}\ot\AAA_{k-i-1})=\AAA_{k,i+1}(\AAA_{i+1,i-1}(\AAA_{i-1}\ot\AAA_{2})\ot\AAA_{k-i-1})=y(\AAA_{i-1}\ot\AAA_2\ot\AAA_{k-i-1})$
with $y=\AAA_{k,i+1}(\AAA_{i+1,i-1}\ot 1)$;  $\AAA_{i+1,i-1}$ acts in the
first $i+1$ positions.
Using  (a) we have 
\begin{align*}
\AAA_k(1\ot T_i)&=y(\AAA_{i-1}\ot\AAA_2\ot \AAA_{k-i-1})(1\ot T_i)
       =y(\AAA_{i-1}\ot\AAA_2(1\ot T)\ot \AAA_{k-i-1})
\\
      &=y(\AAA_{i-1}\ot\AAA_2(T'\ot 1)\ot\AAA_{k-i-1})=\AAA_k(T_i^\prime\ot1).
\end{align*}
(c) Suppose \rf[relianti] is fulfilled for $h$ and $g$ in $\h$. Then it is valid
for $hg$ as well because $\AAA_k(1\ot hg)=\AAA_k(1\ot h)(1\ot g)=\AAA_k(h'\ot
1)(1\ot g)=\AAA_k(h'\ot g)=\AAA_k(1\ot g)(h'\ot 1)=\AAA_k( g' h'\ot
1)=\AAA_k((hg)^\prime\ot 1)$. Since $T_i$ generate $\h$, the proof of  \rf[relianti]   is  complete.
\\
We show that $\AAA_k(z_\lam^\prime\ot z_\lam)$ is essentially
idempotent. Noting that $1\ot z_\lam$ is central in $\H$, using
\rf[relianti] several times, $z_\lam^2=z_\lam$, \rf[d-anti], and
Lemma\,\ref{l-tcentral},  it follows that 
\begin{equation}        \nn
\begin{split}           
(\AAA_k(z_\lam^\prime\ot z_\lam))^2&=\AAA_k(1\ot z_\lam)
\AAA_k(1\ot z_\lam)=\AAA_k^2(1\ot z_\lam)
\\
&=\AAA_k\ts\sum_{w\in\S}(T_w^\prime\ot T_w^\ast z_\lam)
=\AAA_k\ts\sum_{w\in\S}(1\ot T_wT_w^\ast z_\lam)
\\
&=\AAA_k(1\ot t z_\lam)=t_\lam\AAA_k(1\ot z_\lam)
=t_\lam\AAA_k(z_\lam^\prime\ot z_\lam).
\end{split}
\end{equation}
Hence, $\pi_\lam$ is idempotent. Using \rf[relianti] again and $z_\lam
z_\mu=0$ for $\lam\ne\mu$, we have 
\begin{equation}\nn
\pi_\lam\pi_\mu=t_\lam^{-1}t_\mu^{-1}\AAA_k(1\ot z_\lam)
\AAA_k(1\ot z_\mu)=t_\lam^{-1}t_\mu^{-1}\AAA_k^2(1\ot z_\lam z_\mu)=0,
\end{equation}
proving that $\pi_\lam$ and $\pi_\mu$ are orthogonal. Since
$1=\sum_{\lam}z_\lam$, $\sum_\lam\AAA_k(z_\lam^\prime\ot
z_\lam)=\sum_\lam\AAA_k(1\ot z_\lam)=\AAA_k$, and \rf[spec] is proved.
\\
To the last assertion. We briefly write $x$ for the set of commuting variables
$x_1,\dots,x_N$. Changing the role of $w$ and $w^{-1}$ in \rf[e-zlam],
applying then involution~$^\prime$, and using Lemma\,\ref{invzentral}, we have
$z_{\lam'}=\frac{d_\lam}{t_\lam}\sum_{w\in\S}\chi^\lam(T_{w^{-1}})T_w^\prime$.
Applying $\chi^{\lam'}$, dividing by
$\chi^{\lam'}(z_{\lam'})=d_{\lam'}=d_\lam$, multiplying by  $s_\lam(x)s_{\lam'}(x)$, and finally using \rf[spurchar] and \rf[pia], we obtain
\begin{align*}
s_\lam(x)s_{\lam'}(x)
&={\ts\frac{1}{t_\lam}\sum_{w\in\S}}\,\chi^\lam(T_{w^{-1}})\chi^{\lam'}
(T_w^\prime)s_\lam(x)s_{\lam'}(x)
\\
&={\ts\frac{1}{t_\lam}\sum_{w\in\S}}\,
\TR(T_{w^{-1}}z_\lam)\TR(T^\prime_wz_{\lam'})
\\
&=(\TR\ot\TR)\bigl({\ts\frac{1}{t_\lam} \sum_{w\in\S}}(T_w^\prime\ot 
T_{w^{-1}})(z_\lam^\prime \ot z_\lam)\bigr)
\\
&=(\TR\ot\TR)({\ts\frac{1}{t_\lam}}\AAA_kz_\lam^\prime\ot z_\lam)
 =(\TR\ot\TR)\pi_\lam.
\end{align*}
\end{proof}
\begin{remark}
The explicit formula for $t_\lam$ is $t_\lam=k!q^{c(\lam)}\frac{H_q(\lam)
}{h(\lam)}$, where $H_q(\lam)=\QM^{-1}\prod_{x\in\lam}(q^{h(x)}-q^{-h(x)})$ and
$c(\lam)=\sum_{(i,j)\in\lam}(j-i)$, see 
\cite[(1.6)]{HLR}. 
\end{remark}
\section{Corepresentations of $\mathcal{O}(GL_q(N))$ and $\mathcal{O}(SL_q(N))$} \label{s-corep}
In this section we recall notions and facts from the theory of
corepresentations of Hopf algebras. All corepresentations are assumed to be
finite dimensional. For a Hopf algebra $\AA$ the following statements are
equivalent \cite[Theorem\,11.3, p.\,403]{KS}:  (i) Every corepresentation of
$\AA$ is a direct sum of irreducible corepresentations,  (ii) $\AA$ is the
linear span of all matrix elements of all irreducible corepresentations. In
this case we call $\AA$ {\em cosemisimple}. The {\em character}\/ of a matrix
corepresentation $v=(v_j^i)_{i,j=1,\dots,d}$ of $\AA$ is the element
$\chi_v=\sum_{i=1}^dv^i_i$ in $\AA$. In case of the zero corepresentation,
$v=\chi_v=0$. For  $\AA$ cosemisimple two finite dimensional corepresentations
are equivalent if and only if their characters coincide, see
\cite[Corollary\,11.18, p.\,407]{KS}.\\
Two idempotents $p_1$ and $p_2$ of an algebra $A$ are called {\em
  equivalent}\/ if there are elements $a,\,b\in A$ with $ab=p_1$ and
$ba=p_2$. Obviously, two idempotents $p_1$ and $p_2$, $p_1\ne0$ or $p_2\ne0$, are {\em not}\/ equivalent, if $p_1xp_2=0$ for all $x\in A$.
\begin{lemma}          \label{l-character}
Let $\AA$ be a cosemisimple Hopf algebra and  $v=(v^i_j)$ be a matrix 
corepresentation of $\AA$ on the  vector space $V$. 
\\
{\rm (i)} Suppose $P\in\Mor(v)$ is 
idempotent. Then the restriction of $v$ to the image of $P$ defines a
subcorepresentation $v(P)$ of $v$ with character
$\chi_{v(P)}=\sum_{i,j}P^i_jv^j_i$.
\\
{\rm (ii)}  Let $P$ and $Q$ be idempotents in $\Mor(v)$. Then the
corresponding subcorepresentations $v(P)$ and $v(Q)$ are equivalent if and
only if $P$ and $Q$ are equivalent.
\\
{\rm (iii)}  Let $P_i\in\Mor(v)$, $i=0,\dots,m$, $P_0\ne0$, be 
equivalent idempotents and let $P=\sum_i\alpha_iP_i$, $\alpha_i\in\C$, be
idempotent. Then there exist non-negative integers $r$ and $s$, $s\ne0$, with
$\sum_i \alpha_i=s^{-1}r$. Moreover, we have the following equivalence of
subcorepresentations of $v$: $s\cdot v(P)\cong r\cdot v(P_0)$.
\end{lemma}
\begin{proof} Throughout the  proof we sum over repeated indices. (i) We
determine the matrix coefficients of $w=v(P)$. Let $\{e_i\}$ and $\{f_J\}$
be bases of $V$ and $W{=}\im P$, resp.\  Define linear mappings  $A\colon
W\to V$, $Aw=w$, and $B\colon V\to W$, $Bv=Pv$. Obviously, $P=AB$ and
$BA=\id_W$. The matrix elements of $A$ and $B$ corresponding to the chosen
bases are determined by $f_J=\sum A^j_Je_j$ and $B(e_i)=\sum B^K_if_K$,
resp.\  Let $\vp$ denote the associated to $v$ right comodule mapping ,
$\vp(e_i)=\sum e_j\ot v^j_i$. Since $A=PA$, $P\in\Mor(v)$, and $P=AB$ we obtain 
\begin{equation}\nn
\begin{split} 
\vp(f_N)&={\ts \sum}\,\vp(A^j_Ne_j)={\ts \sum}\,e_i\ot v^i_jA_N^j=
{\ts \sum}\,e_i\ot v^i_jP^j_kA_N^k
\\
&={\ts \sum}\,e_i\ot v^j_kP^i_jA_N^k={\ts \sum}\,A^i_Me_i\ot
B^M_jv^j_kA^k_N={\ts\sum}\,f_M\ot(BvA)^M_N.
\end{split}
\end{equation}
Hence, $W=\im P$ defines a subcorepresentation with matrix elements
$w=BvA$. Therefore,  $\chi_w=\sum B^M_jv^j_kA^k_M=\sum (AB)^k_jv^j_k=\sum
P^k_jv^j_k$. 
\\
(ii) Suppose $P\sim Q$ are equivalent idempotents, i.\,e.\  $P=AB$ and $Q=BA$
for some $A,\,B\in \Mor(v)$. By (i) the characters of the subcorepresentations
$v(P)$ and $v(Q)$ are $\chi_{v(P)}=\sum_{k,l}P^k_lv^l_k$ and $\chi_{v(Q)}=\sum
Q^k_lv^l_k$, resp. Inserting $P=AB$, $Q=BA$, and  using $B\in\Mor(v)$ gives
$\chi_{v(P)}=\sum A^k_xB^x_lv^l_k=\sum A^k_xv^x_lB^l_k=\sum
Q^l_xv^x_l=\chi_{v(Q)}$. Since $\AA$ is cosemisimple, $v(P)\cong
v(Q)$. Suppose now that $v(P)\cong v(Q)$ are equivalent corepresentations,
i.\,e.\  there is a bijective map $J$, $J\in\Mor(v(P),v(Q))$. Let $W_1$ and
$W_2$ denote the images of $P$  and $Q$, resp.  Moreover, let $A_1,B_1$ and
$A_2,B_2$ be the corresponding mappings from the proof of (i) for $P$ and $Q$,
resp. Since $J\in\Mor(v(P),v(Q))$, $J(B_1vA_1)=(B_2vA_2)J$. Applying the
counit $\ve$ and choosing $A=A_1J^{-1}B_2$ and $B=A_2JB_1$ one gets $AB=P$ and
$BA=Q$. Hence, $P\sim Q$.
\\
(iii) By (i) and (ii):
\begin{equation}  \label{char}
\chi_{v(P)}=
\ts\sum P^k_lv^l_k=\ts\sum\alpha_i(P_i)^k_lv^l_k=
\ts\sum\alpha_i\chi_{v(P_i)}=
(\ts\sum_i\alpha_i)\chi_{v(P_0)}.
\end{equation} 
Since $\AA$ is cosemisimple there exist  integers $l$ and $s_j$, $s_j\ge1$,
$j=1,\dots,l$, and  irreducible corepresentations $\vp_j\ne0$, $j=1,\dots,l$,
such that $v(P_0)\cong\sum_j s_j\vp_j$. Set $r_j=s_j\sum_i\alpha_i$,
$j=1,\dots,l$. {}From  \rf[char] it follows
$\chi_{v(P)}=\sum_jr_j\chi_{\vp_j}$. Since $\AA$ is cosemisimple and since
irreducible characters are linearly independent, $r_j$ are non-negative
integers and $v(P)\cong\sum_jr_j\vp_j$. Hence $\sum_i\alpha_i=s^{-1}_1r_1$ is
rational and $s_1v(P)\cong r_1 v(P_0)$. 
\end{proof}
Now let $\AA$ be one of the Hopf algebras $\oglqn$ or  $\oslqn$. Since  $q$ is
transcendental, $\AA$ is cosemisimple, see \cite[Theorem\,11.22,
p.\,410]{KS}. Throughout let $\vp$ and $\psi$ denote the corepresentations
$u^{\ot k}$ and $(u^\fettc)^{\ot k}$, respectively. By \rf[e-mor1] below it is
obvious that $\vr(\h)\subseteq\Mor(\vp)$ and
$\vr_\fettc(\h)\subseteq\Mor(\psi)$. Let $p\in\h$ be  idempotent, $P=\vr(p)$,
and $P^\fettc=\vr_\fettc(p)$. Since, $P\in\Mor(\vp)$ (resp.\
$P^\fettc\in\Mor(\psi)$) and by Lemma\,\ref{l-character}\,(i), the restriction
of the corepresentation $\vp$  (resp.\  of $\psi$) to the image of $P$ (resp.\
of $P^\fettc$) defines a subcorepresentation  $\vp(P)$ (resp.\
$\psi(P^\fettc)$), perhaps $0$. In particular, let $p_\lam$ and $\ov{p}_\lam$
be  minimal subidempotents of $z_\lam$. Since $z_\lam\h$ is a simple ideal,
$p_\lam\sim\ov{p}_\lam$ are equivalent idempotents. Thus
$\vr(p_\lam)\sim\vr(\ov{p}_\lam)$ (resp.\
$\vr_\fettc(p_\lam)\sim\vr_\fettc(\ov{p}_\lam)$) are equivalent idempotents
too. Hence by  Lemma\,\ref{l-character}\,(ii),
$\vp(\vr(p_\lam))\cong\vp(\vr(\ov{p}_\lam))$ (resp.\
$\psi(\vr_\fettc(p_\lam))\cong\psi(\vr_\fettc(\ov{p}_\lam))$) are equivalent
corepresentations.   Let  $\vp_\lam$ (resp.\ $\psi_\lam$) denote this equivalence class.
We determine the dimensions  of $\vp_\lam$ and  $\psi_\lam$. Let $\TR_1$
denote the evaluation of $\TR$ at $x_1=\dots=x_N=1$, see  \rf[d-tau]. By
Remark\,\ref{hecke} and by $\sum_{c\models k}E_c=I$ we have 
\begin{equation}\label{tr1}
\TR_1=\tr{\circ}\vr=\tr{\circ}\vr_\fettc.
\end{equation}
By \rf[spurpl],
$\rank(\vr(p_\lam^i))=\rank(\vr_\fettc(p_\lam^i))=\TR_1(p_\lam^i)=s_\lam(1,\dots,1)$.
The explicit value is 
\begin{equation}                     \label{dimlam}
\delta_\lam(N):=\dim\vp_\lam=\dim\psi_\lam=s_\lam(1,\dots,1)=
h(\lam)^{-1}\prod_{(i,j)\in\lam}(N+j-i),
\end{equation}
see \cite[I.3~Example~4, p.~45]{Mac}.
\begin{lemma}       \label{l-contra}
Let $\AA$ be one of the Hopf algebras  $\oglqn$ or  $\oslqn$. Assume $p\in\h$
be idempotent, $P^\ast=\vr(p^\ast)$, and $P^\fettc=\vr_\fettc(p)$. Then  the
following corepresentations  are equivalent
\begin{equation}
\vp(P^\ast)^\fettc\cong \psi(P^\fettc).
\end{equation}
\end{lemma}
\begin{proof}
We first prove
\begin{equation}           \label{rcrc}
\vr_\fettc(h)^{\iinv}_{\jinv}=\vr(T_{\wo}h^\ast T_\wo^{-1}) ^{\jvec}_{\ivec},
\end{equation}
for  $h\in\h$, $\ivec=(i_1,\dots,i_k)$, $\jvec=(j_1,\dots,j_k)$,
$\iinv=(i_k,\dots,i_1)$, and $\jinv=(j_k,\dots,j_1)$. First let $h=T_n$ be a
single generator. By \rf[bwo], $b_nb_\wo=b_\wo b_{k-n}$. Applying the
representation $b_n\mapsto T_n$ to this equation and recalling $T_n^\ast=T_n$,
we have
\begin{equation}  \nn
\vr(T_\wo T_n^\ast T_\wo^{-1})^{\jvec}_{\ivec}=\vr(T_{k-n})^{\jvec}_{\ivec} 
=\Rda[i_{k-n},i_{k-n+1}]{j_{k-n},j_{k-n+1}}
=\Rch[j_{k-n+1},j_{k-n}]{i_{k-n+1},i_{k-n}}=\vr_\fettc(T_n)^\iinv_\jinv;
\end{equation}
that proves \rf[rcrc] for $h=T_n$. Suppose now \rf[rcrc] is valid for
$g,h\in\h$. By representation property of $\vr$ and $\vr_\fettc$ and
antimultiplicativity of $^\ast$, we finally conclude
\begin{equation}\nn
\vr_\fettc(gh)^\iinv_\jinv ={\ts\sum_\xvec}\vr_\fettc(g)^\iinv_\xinv
\vr_\fettc(h)^\xinv_\jinv={\ts\sum_\xvec}
\vr(T_\wo g^\ast T_\wo^{-1}) _\ivec^\xvec 
\vr(T_\wo h^\ast T_\wo^{-1})_\xvec^\jvec=
\vr(T_\wo (gh)^\ast T_\wo^{-1})^\jvec_\ivec.
\end{equation}
Since $\{T_n\}$ generate $\h$, \rf[rcrc] is proved. We sum over repeated
multi-indices. By Lemma\,\ref{l-character}\,(i), \rf[rcrc], and
$S(\chi_v)=\chi_{v^\fettc}$ we have
\begin{equation}  \nn
\begin{split}
\chi_{\psi(P^\fettc)}&={\ts\sum}(P^\fettc)^\ivec_\jvec\,(u^\fettc)
^\jvec_\ivec =\ts\sum\vr(T_\wo p^\ast T_\wo^{-1})^\jinv_\iinv(u^\fettc)^\jvec_\ivec
\\
&=\ts\sum\vr(T_\wo)^\jinv_\xvec (P^\ast)^\xvec_\zvec\vr(T_\wo^{-1})_\iinv 
^\zvec S(u^\iinv_\jinv)
=S\bigl(\ts\sum \vr(T_\wo)^\jinv_\xvec (P^\ast)^\xvec_\zvec \vr(T_\wo^{-1})^\iinv_\jinv \,u^\zvec_\iinv\bigr)
\\
&=
S(\ts\sum (P^\ast)^\iinv_\zvec u^\zvec_\iinv) =S(\chi_{\vp(P^\ast)}) =
\chi_{\vp(P^\ast)^\fettc}.
\end{split}
\end{equation}
Since $\AA$ is cosemisimple $\vp(P^\ast)^\fettc$ and $\psi(P^\fettc)$ are
equivalent corepresentations. 
\end{proof}
\begin{cor}\label{c-contra}
$\vp^\fettc_\lam\cong\psi_\lam$.
\end{cor}
\begin{proof}
Let $p_\lam$ be a minimal subidempotent of $z_\lam$. Since
$z_\lam=z_\lam^\ast$, $p_\lam\sim p_\lam^\ast$. Hence
$\vp_\lam^\fettc=\vp(\vr(p_\lam))^\fettc\cong\vp(\vr(p_\lam^\ast))^\fettc\cong\psi(\vr_\fettc(p_\lam))=\psi_\lam$.
\end{proof}
We recall the  Brauer-Schur-Weyl duality for quantum groups  of type A,
cf.~\cite{Ha} or \cite[Theorem~8.38, Proposition~11.20, and
Proposition~11.21]{KS}.
\begin{prp}
\label{bsw}
Let $\AA$ be one of the Hopf algebras $\oglqn$ or $\oslqn$. Let $q$ be a
transcendental complex number, $k\in \N$, and $\lam\vdash k$.
\\
{\rm (i)} The representation $\vr\colon\h\to \Mor(u^{\ot k})$ is surjective. 
The subcorepresentation   $\vp_\lam$ of $\vp$ is irreducible.
\\
{\rm (ii)} $\ker\vr=\bigoplus\limits_{\lam\vdash
    k,\lam_1^\prime>N}z_\lam\h$. 
The subcorepresentation $\vp_\lam$ is zero if and only if $\lam_1^\prime>N$.
\end{prp}
\begin{cor}\label{c-bsw}
Let $\AA$ and $q$ be as above, and let $\lam,\mu\vdash k$ be partitions of
$k$. Then $\vp_\lam\cong\vp_\mu$ if and only if $\lam=\mu$ or
both $\lam_1^\prime>N$  and $\mu^\prime_1>N$.
\end{cor}
\begin{proof}
$\leftarrow$ is trivial. $\to$: Suppose $\lam_1^\prime\le N$ or
$\mu_1^\prime\le N$ and $\lam\ne\mu$. Let $p_\lam$ and $p_\mu$ denote minimal
idempotents corresponding to $\lam$ and $\mu$, respectively, and 
$P_\lam=\vr(p_\lam)$, $P_\mu=\vr(p_\mu)$. Since
$P_\lam\ne0$ or $P_\mu\ne0$ by Proposition\,\ref{bsw}\,(ii) and since $P_\lam
X P_\mu=0$, $X\in\Mor(u^{\ot k})$, $P_\lam$ and $P_\mu$
are inequivalent.  By Lemma~\ref{l-character}~(ii), $\vp_\lam$ and $\vp_\mu$
are inequivalent too.
\end{proof}
Corresponding to $\Gamm=\gtz$ consider the representation  $\vr_\tau$  of $\H$
on $V^{\ot 2k}$, $\vr_+=\vr_\fettc\ott\vr$ and
$\vr_-=(\vr_\fettc\ott\vr)\alpha$, where $\alpha:\H\to\H$ denotes the
automorphism $\alpha(h\ot g)=g\ot h$, $g,h\in\h$. The next Proposition settles
the problem of determining the rank of the antisymmetrizer
$a_k=\vr_\tau(\AAA_k)$. The basic result is
\begin{prp}
\label{p-dec}
The idempotent $\varPi_\lam^\tau=\vr_\tau(\pi_\lam)$ defines a
subcorepresentation $\phi(\varPi_\lam^\tau)$ of $\phi=(u^\fettc)^{\ot k} \ot
u^{\ot k}$ with 
\begin{equation}\label{e-corep}
\phi(\varPi_\lam^+)\cong\psi_{\lam^\prime}\ot \vp_\lam
\quad\text{and}\quad
\phi(\varPi_\lam^-)\cong\psi_\lam\ot \vp_{\lam^\prime}.
\end{equation}
\end{prp}
\begin{proof}
(a) We first calculate the rank of $\varPi_\lam^\tau$ in $V^{\ot 2k}$, i.\,e.\
the dimension of the corepresentation $\phi(\varPi_\lam^\tau)$. Let $\pi\in\H$
be idempotent. By \rf[tr1], $\rank\vr_\tau(\pi)=(\TR_1\ot\TR_1)\pi$. In
particular for $\pi=\pi_\lam$  by \rf[ttPi] and  \rf[dimlam] we have
\begin{equation}\label{rank}
\rank\varPi_\lam^\tau =\delta_{\lam^\prime}(N)\delta_\lam(N).
\end{equation}
(b) Let $\{p_\lam^{ij}\colon i,j=1,\dots,d_\lam\}$ be a linear  basis of
$z_\lam\h$ consisting of mutually equivalent idempotents. Then
$\{(p_\lam^{ij})^\prime\colon i,j=1,\dots,d_\lam\}$ is a basis of
$z_{\lam'}\h$. By \rf[pia], $\pi_\lam$ is a subidempotent of
$z_{\lam^\prime}\ot z_\lam$. Hence there are complex numbers
$\alpha_{rsmn}^\lam$, $r,s,m,n=1,\dots,d_\lam$, such that
$\pi_\lam=\sum\alpha^\lam_{rsmn}(p_{\lam}^{rs})^\prime\ot
p^{mn}_\lam$. Moreover, the idempotents $(p_{\lam}^{rs})^\prime\ot
p^{mn}_\lam$ are mutually equivalent  in $\H$. Applying $\vr_+$ (resp.\
$\vr_-$) gives $\varPi_\lam^+=\sum
\alpha^\lam_{rsmn}\vr_\fettc((p_{\lam}^{rs})^\prime)\ot \vr(p^{mn}_\lam)$ (resp.\
$\varPi_\lam^-=\sum\alpha^\lam_{rsmn}\vr_\fettc(p_{\lam}^{mn})\ot
\vr((p^{rs}_\lam)^\prime)$). One verifies that $v=\phi$, $P=\varPi_\lam^+$ (resp.\
$P=\varPi_\lam^-$), and $P_{i}=\vr_\fettc((p_{\lam}^{rs})^\prime)\ot
\vr(p^{mn}_\lam)$ (resp.\
$P_i=\vr_\fettc(p_{\lam}^{mn})\ot\vr((p^{rs}_\lam)^\prime)$) satisfy the
assumptions of Lemma\,\ref{l-character}\,(iii). Hence,
$\phi(\varPi_\lam^+)\cong\alpha^\lam\, \psi_{\lam'}\ot\vp_\lam$, (resp.\
$\phi(\varPi_\lam^-)\cong\alpha^\lam\, \psi_{\lam}\ot\vp_{\lam'}$), where
$\alpha^\lam=\sum\alpha^\lam_{rsmn}$ is rational. In particular, by
\rf[dimlam],
$\dim\phi(\varPi_\lam^\tau)=\alpha^\lam\,\delta_{\lam'}(N)\delta_\lam(N)$.
Comparing  this  with  \rf[rank] gives $\alpha^\lam=1$; the proof is complete. 
\end{proof} 
The main step to simplify the study of the algebras $\gdwl$ and $\gdwi$ is the
reduction of the antisymmetrizer $\ant$ of $\gol[k]$ to the antisymmetrizer
$a_k$ of $\vr_\tau(\H)$. We provide an isomorphism of the right
$\AA$-comodules $\uckuk$ and $\dr{\uhr}\gol[k]$   that maps $a_k$ into
$\ant$. Here and in the remainder of the article $V$ always denotes the
$N$-dimensional complex vector space $\C^N$ with canonical basis
$e_1,\dots,e_N$.
\\
Set $\fr[+]{s}=1+\fr{s}+q^{-2N-2}$, $\fr[-]{s}=\fr{s}-q^{-2}-q^{-2N}$, and
define the $N^2{\times}N^2$-matrices $(\grave{R}^\pm)^{ab}_{rs}=q^{2s-2a}
(\Rda{\pm1})^{bs}_{ar}$. We recall some well-known properties of the  matrices
$\Rre$ and $\grave{R}^\pm$, see \cite[Lemma\,3.3 and Lemma\,3.4]{SS2}:
\begin{alignat}{3} 
\Mor(u\ott u)&=\langle\Rda{},\Rda{-1}\rangle,&
\Mor(u^\fettc\ott u^\fettc)&=\langle\Rch{+}, \Rch{-}\rangle, &
\label{e-mor1}\\
  \Mor(u^\fettc\ott u,u\ott  u^\fettc)&=\langle\Rre^{+},\Rre^{-}\rangle,&
  \Mor(u\ott u^\fettc,u^\fettc\ott u)&=\langle \Rli{+},\Rli{-} \rangle, &
\nn\\ 
\Rre^\pm\grave{R}^\mp&=\grave{R}^\mp\Rre^\pm=I,&
 \sum_{i,j}q^{2s-2j}(P_\pm)^{in}_{jm}(P_\pm)^{js}_{ir}&=
\QP^{-2}(\delta_{nr}\delta_{ms}+\fr[\pm]{s}q^{2s}\delta_{nm}\delta_{rs}),&                       \label{e-mor2} \\ 
\Rre[23]^\pm\Rre[12]^\pm \Rda[23]{}&=\Rda[12]{}\Rre[23]^\pm\Rre[12]^\pm,&
\Rre[12]^\pm\Rre[23]^\pm \Rch[12]{}&=\Rch[23]{}\Rre[12]^\pm\Rre[23]^\pm.& 
\label{e-mor3} 
\end{alignat} 
\begin{prp}                                   \label{p-mor}
Let $\AA$ denote one of the Hopf algebras $\oglqn$ or $\oslqn$, $\Gamm=\gtz$,
$\phi=\uckuk$, and $k\ge1$.
\\ {\rm (i)}
There exists an isomorphisms $\ip{k}$ of right comodules $\phi$ and
$\dr{\uhr}\gol[k]$ such that
\begin{equation}               \label{e-antianti} 
\ip{k}a_k=\ant\ip{k}. 
\end{equation} 
{\rm(ii)} 
There exists an isomorphism  $\idp{k}$ of right comodules $\phi{\uhr}\im  a_k$
and $\dr{\uhr}\gdwl[k]$.
\\ {\rm (iii)} There exists an isomorphism  of $\gdwi[k]$ and the vector space
$(\im a_k)^\phi$ of  $\phi$-invariant elements of $\im a_k$. 
\end{prp} 
\begin{proof} 
To simplify notations, throughout the proof we shall write $C_n$ instead of
$C_{n,n+1}$ for $C\in\End(V\ot V)$ acting in position $n$ and $n+1$ of the
tensor product $V^{\ot k}$, $n\le k-1$.
\\ 
(i) Let $\{e_\ivec\colon \ivec=(i_1,\dots,i_{2k})\}$ and
$\{\theta_{i_1i_2}\ota\cdots\ota\theta_{i_{2k-1}i_{2k}}\}$ be the canonical bases of
$V^{\ot 2k}$ and $\gol[k]$, resp. With respect to these bases define the
matrix $\ip{k}=\ipm{k}$ as follows. Set $\ipm{1}=I$ and for $k\ge2$ let
\begin{equation}\label{e-iso}
\ipm[2\,\cdots \,2k-1]{k}= 
\Rrepm[2k-2](\Rrepm[2k-4]\Rrepm[2k-3])\cdots 
(\Rrepm[2i]\cdots\Rrepm[k+i-1])\cdots(\Rrepm[2]\cdots\Rrepm[k]). 
\end{equation} 
The index  $(2,\dots,2k-1)$ indicates that $\ipm{k}$ effectively acts at
positions $2,\dots, 2k-1$ and leaves the first and last coordinates
unchanged. Since $\Rrepm$ is invertible by \rf[e-mor2], $\ipm{k}$ is
invertible. The following two recursion formulas are easily checked
\begin{align*}
\ipm{k}&= 
\ipm[4\,\cdots\,2k-1]{k-1}\Rrepm[2]\cdots\Rrepm[k], \\ 
\ipm{k}&=\ipm[2\,\cdots\,2k-3]{k-1}\Rrepm[2k-2]\cdots\Rrepm[k]. 
\end{align*}
Using  \rf[e-mor1]  and one of the above recursion formulas one shows by
induction on $k$ that the matrix
$\ipm{k}$ defines an isomorphism of corepresentations $\uckuk$ and
$\ucuk$. Since $\dr{\uhr}\gol[k]\cong\ucuk$ by definition of $\dr$ and
\rf[right], the first part of  (i) is proved. Next we shall show equation
\rf[e-antianti]. Since both antisymmetrizer $ a_k$ and $\ant$ are homomorphic
images of the braid group antisymmetrizer $\AAAA_k$ (under
$\vr_\tau\circ\bsig$ and $\gamma$, resp.) it suffices to prove
$\ipm{k}\vr_\tau(T_n^{-1}\ot T_n)=(\sig_\tau)_n\ipm{k}$ for
$n=1,\dots,k-1$. By \rf[sig] the equivalent matrix notation of this identity
is
\begin{equation}\label{sign} 
\ipm{k}\Rch[n]{-\tau}\Rda[k+n]{\tau}=\Rrepm[2n]\Rch[2n-1]{-\tau} 
\Rda[2n+1]{\tau}(\Rrepm[2n])^{-1}\ipm{k}. 
\end{equation}
We will  prove \rf[sign] for $\ip[+]{k}$ by induction on $k$. The proof for
$\ip[-]{k}$ is analogous. Since $\ip[+]{2}=\Rre[2]{}$, \rf[sign] is obvious
for $k=2$ and $n=1$.
\\
{\em Case 1}. $n=1$. By the second recursion equation, \rf[e-mor3], $k\ge3$, and induction
assumption  we have
\begin{equation*}
\begin{split}
 \ip[+]{k}\Rchm[1]\Rda[k+1]{} 
&=\ip[+\,2\,\cdots\,2k-3]{k-1}\Rre[2k-2]\cdots \Rre[k+2]\Rre[k+1]\Rre[k]\Rchm[1] 
\Rda[k+1]{}
 \\ &=\ip[+\,2\,\cdots\,2k-3]{k-1}\Rchm[1]\Rre[2k-2]\cdots\Rre[k+2] 
\Rda[k]{}\Rre[k+1]\Rre[k]
 \\ &=\ip[+\,2\,\cdots\,2k-3]{k-1}\Rchm[1] 
\Rda[k]{}\Rre[2k-2]\cdots \Rre[k] 
\\ 
&=\Rre[2]\Rchm[1]\Rda[3]{}{\Rre[2]}^{-1}
\ip[+\,2\,\cdots\,2k-3]{k-1} \Rre[2k-2]\cdots 
\Rre[k]
 \\ &=\Rre[2]\Rchm[1]\Rda[3]{}{\Rre[2]}^{-1}\ip[+]{k}. \end{split} 
\end{equation*}
{\em Case 2}. $2\le n\le k-1$. By the first recursion equation, \rf[e-mor3], $k\ge3$, and induction
assumption we have
\begin{equation}\nn
\begin{split} 
\ip[+]{k}\Rchm[n]\Rda[k+n]{} &= 
\ip[+\,4\,\cdots\,2k-1]{k-1}\Rre[2]\cdots\Rre[k]\Rchm[n]\Rda[k+n]{} 
\\ &= 
\ip[+\,4\,\cdots\,2k-1]{k-1}\Rda[k+n]{}\Rre[n]\cdots\Rre[n]\Rre[n+1]\Rchm[n] 
\Rre[n+2]\cdots\Rre[k] 
\\ &= \ip[+\,4\,\cdots\,2k-1]{k-1} 
\Rchm[n+1]\Rda[k-1+n+1]{}\Rre[2]\cdots\Rre[k] \\ 
&=\Rre[2n]\Rchm[2n-1]\Rda[2n+1]{}\Rre[2n] \ip[+\,4\,\cdots\,2k-1]{k-1} 
\Rre[2]\cdots\Rre[k] \\ &=\Rre[2n]\Rchm[2n-1]\Rda[2n+1]{}\Rre[2n]\ip[+]{k}. 
\end{split} 
\end{equation} 
The proof of (i) is complete.
\\
(ii) By (i), $\ip{k}\colon V^{\ot 2k}\to\gol[k]$ is a linear
isomorphism. Moreover, the kernels of $ a_k$ and $\ant$ are bijectively mapped
into each other. Consequently, $\ip{k}$ can be factorised to an isomorphism
$V^{\ot 2k}{/}\ker a_k\to\gol[k]{/}\ker\ant$. Since $\im a_k$ is also
$\phi$\,-covariant and since $\im a_k\cong V^{\ot 2k}{/}\ker a_k$,
$\idp{k}\colon\im a_k\to\gdwl[k]$ is again an isomorphism of right comodules.
\\
(iii)  Since $\AA$ is cosemisimple and $\idp{k}$ is an isomorphism of right
comodules, $\idp{k}$ is the direct sum of isomorphisms of the isotypical
components. In particular, the trivial components (corresponding to the
trivial corepresentation $1$) are bijectively mapped into each other. Hence, the
restriction of $\idp{k}$ to the space of invariant elements of $\im a_k$ is an
isomorphism to $\gdwi[k]$. 
\end{proof} 
\subsection*{Proofs of Theorem~\ref{left} and 
Theorem\,\ref{bi}\,(i)}
We first show that  the following combinatorial formula holds:
\begin{equation}\label{e-dim}
 \sum_{\lam\vdash k}\delta_{\lam'}(N)\delta_\lam(N)=\binom{N^2}{k},
\end{equation} 
where $N$ and $k$ are non-negative integers. Use formula
\cite[$(4.3')$, I.\,4, p.\,65]{Mac}: $\prod(1+x_iy_j)=\sum_\lam
s_\lam(x)s_{\lam'}(y)$. We restrict it to the finite set of variables $x_1,\dots,x_N,\,y_1,\dots, y_N$
and consider the natural bi-grading of  polynomials in these
variables. Comparing the homogeneous components of degree $(k,k)$, inserting
$x_1=\dots=y_N=1$, and using \rf[dimlam] gives \rf[e-dim]. By \rf[rank] we
finally obtain
\begin{equation}\nn
\dim\gdwl[k]=\sum_{\lam\vdash k}\rank \varPi^\tau_\lam=\sum_{\lam\vdash k}
\delta_{\lam'}(N)\delta_\lam(N)=\binom{N^2}{k}.
\end{equation}
In particular,  the only partition $\lam$ of $N^2$ with
$\rank\varPi^\tau_\lam\ne0$ is $\lam=(N^N)$. In this case $\lam=\lam'$ is
symmetric. By the above formula, $\dim\gdwl[N^2]=1$. Hence, there exists a
unique up to scalars form of maximal degree $N^2$; the proof of
Theorem\,\ref{left} is complete.
\\
To prove Theorem\,\ref{bi}\,(i) consider the formal power series
$p(t)=(1+t)(1+t^3)(1+t^5)\cdots=1+c_1t+c_2t^2+\cdots$. Then
$c_k=\Card\{\lam:\lam\vdash k,\lam=\lam'\}$. In fact, let $\lam\vdash k$,
$\lam=\lam'$, be a symmetric Young diagram. Denote the hook consisting of the
first row and first column of $\lam$ by $\lam^{(1)}$ and the remaining part of
$\lam$ by $\lam_{(1)}$. Obviously, $\lam_{(1)}$ is again a symmetric
diagram. Repeating the above procedure we step by step get a sequence of
symmetric hooks $\lam^{(1)},\dots,\lam^{(r)}$ with $k=\sum_in_i$,
$n_i=|\lam^{(i)}|$, $n_1>\cdots>n_r$, and $n_i$ is odd for each $i$. Now we
associate to $\lam$ the expression $t^{n_1}\cdots t^{n_r}$ appearing in
$p(t)$. Conversely, to each summand $t^{n_1}\cdots t^{n_r}$ in $p(t)$ we
associate a symmetric Young diagram by putting together  symmetric hooks of
weights $n_i$. Similarly, the polynomial $s(t)=(1+t)(1+t^3)\cdots
(1+t^{2N-1})=1+c_1t+\cdots +c_{N^2}t^{N^2}$  has coefficients
$c_k=\Card\{\lam:\lam\vdash k,\lam=\lam',\lam_1^\prime\le N\}$.
\\
Since $\im\varPi_\lam^\tau=0$ for $\lam_1^\prime>N$ by \rf[rank] and
\rf[dimlam], it suffices to consider the elements of
$\PP_k:=\{\lam\colon\lam\vdash k,\,\lam^\prime_1\le N\}$. Then $\vp_\lam$ is
nonzero by Proposition~\ref{bsw}~(ii). We collect three simple facts from the
theory of corepresentations: Let $W$ be the vector space where the
corepresentation $w$ acts. Then the vector space $W^w$ of under $w$ invariant
elements  is isomorphic to $\Mor(1,w)$. Further, $\Mor(1,v\ot
w)\cong\Mor(v^\fettc,w)$, and $v\cong v^{\fettc\fettc}$ (since  $\AA$ is
cosemisimple).
\\
By Proposition\,\ref{p-mor}\,(iii),  \rf[e-corep], and
Corollary\,\ref{c-contra} we have for $\Gamm=\Gamm_{+,z}$
\begin{align*}
\dim(\gdwi[k])&=\sum_{\lam\in\PP_k}\dim(\im\varPi^+_\lam)^\phi=\sum_{\lam\in\PP_k}
\dim\Mor(1,\phi(\varPi^+_\lam))
\\
&=\sum_{\lam\in\PP_k}\dim\Mor(1,\psi_{\lam'}\ot\vp_\lam)=
\sum_{\lam\in\PP_k}\dim\Mor(\psi_{\lam'}^\fettc,\vp_\lam)
\\
&=\sum_{\lam\in\PP_k}\dim\Mor(\vp_{\lam'}^{\fettc\fettc},\vp_\lam)
=\sum_{\lam\in\PP_k}\dim\Mor(\vp_{\lam'},\vp_\lam)
\\
&=\Card\{\lam:\lam\in\PP_k,\vp_\lam\cong\vp_{\lam'}\}=
\Card\{\lam:\lam\in\PP_k,\lam=\lam'\}.
\end{align*}
The last two steps are  by Schur's Lemma and by  Corollary~\ref{c-bsw}. This
is exactly the coefficient $c_k$ of the polynomial $s(t)$. The proof for the
calculi $\Gamm=\Gamm_{-,z}$ is  analogous; this finishes the proof of
Theorem\,\ref{bi}\,(i).
\subsection*{Proofs of Theorem\,\ref{bi}~(ii)--(iv)}
We begin with  a rather general relation between the antipode of
$(\gten)_\ii$ and the action of the longest word $\sig_\wo$ on bi-invariant
elements. Then we determine this action for $\Gamm=\gtz$. Let $\{\theta_i\}$
and $\{\eta_i\}$ be bases of the linear spaces $\gl$ and $\gr$, respectively,
which are related via $\theta_i=\sum_j\eta_jv^j_i$ and
$\eta_i=\sum_j\theta_jS(v^j_i)$, see \cite[(2.38)]{Wo2} 
\begin{prp}          
\label{antipode} Let $\AA$ be a Hopf algebra and let $\Gamm$ 
be a  bicovariant bimodule over $\AA$ with braiding $\sig$. Let $S$ be the
antipode of the $\NO$-graded super Hopf algebra $\gten$ and let
$\rho\in\Gamm^{\ot k}_\ii$ be  bi-invariant. Then 
\begin{equation}    \label{anti-e} 
S(\rho)=(-1)^{\half k(k+1)}\sig_\wo(\rho). 
\end{equation} 
This equation is also valid for $\rho\in\gdwi[k]$.
\end{prp} 
\begin{proof} Throughout the proof we sum over repeated indices. (a)  Since
$S(\theta)=-\sum_{(\theta)}S(\theta_{(-1)})\theta_{(0)}S(\theta_{(1)})$ and
$\dr(\theta_i)=\theta_j\ot v^j_i$, $S(\theta_i)=-\eta_i$. Using \rf[srr] we
get $S(\theta_\ivec)=(-1)^ {\half k(k+1)}\eta_{\iinv}$. Here
$\iinv=(i_k,\dots,i_1)$ denotes the reversed index sequence of
$\ivec=(i_1,\dots,i_k)$,
$\theta_\ivec=\theta_{i_1}\ota\cdots\ota\theta_{i_k}$, and
$\eta_\jinv=\eta_{j_k}\ota\cdots\ota\eta_{j_1}$.
\\ 
(b) We use induction over $k$,  $k\ge2$, to show 
\begin{equation}        \label{eta-v} \eta_\jinv\, 
v^{\jvec}_{\ivec}=\sig_\wo(\theta_\ivec). 
\end{equation} 
By $\eta_jv^j_i=\theta_i$, $\theta_i\rac a=f^i_n(a)\theta_n$, and \rf[fv], we
have  $\eta_{j_2}\ota\eta_{j_1}v_{i_1}^{j_1}v_{i_2}^{j_2}=\eta_{j_2}\ota
\theta_{i_1}\,v^{j_2}_{i_2}=\eta_{j_2}v_x^{j_2}f_y^{i_1}(v_{i_2}^x)\ota
\theta_y=\sig_{i_1i_2}^{xy}\theta_x\ota\theta_y=\sig(\theta_{i_1}\ota\theta_{i_2})$;
this proves  \rf[eta-v] in case $k=2$. Suppose now \rf[eta-v] holds for
$k$. We prove it for  $k+1$. Let $\wokk$ be the longest word in  $\S[k+1]$,
$\wokk=s_1s_2\cdots  s_k\wo$. Since  $\ell(\wokk)=k+\ell(\wo)$ we  can lift
this equation to braids; afterwards we can apply
$\gamma$. Hence, $\sig_\wokk=\sig_1\cdots\sig_k\sig_\wo$. Making repeatedly use
of $\theta_iv^n_j=\sig_{ij}^{xy}v^n_x\theta_y$ gives
\begin{equation}     \label{nu-v} \theta_\jvec\, 
v_i^j=(\sig_1\cdots\sig_k)^{n\, \nvec}_{\jvec\, i}v^j_n\theta_\nvec. 
\end{equation} 
Using induction assumption, \rf[nu-v], and the above formula for $\sig_\wokk$,
one has
\begin{alignat}{2} \eta_{j}\ota\eta_{\jinv}v^\jvec_\ivec\, v^j_i =&\, 
\eta_j(\sig_\wo)_\ivec^\jvec\,\ota\,\theta_\jvec\, v_i^j &=&\, 
\eta_j(\sig_\wo)_\ivec^\jvec(\sig_1\cdots\sig_k)_{\jvec\, i}^ {n 
\nvec}v^j_n\,\ota\,\theta_{\nvec} \nn \\ =&\,\eta_jv^j_n(\sig_\wokk)^{n 
\nvec}_{\ivec\, i}\,\ota\,\theta_\nvec
&=&\,(\sig_\wokk)^{n\,\nvec}_{\ivec\,i}\theta_n\,\ota \theta_{\nvec}= \sig_\wokk(\theta_{\ivec}\,\ota\theta_{i})\nn 
\end{alignat} 
and the proof of \rf[eta-v] is complete.
\\(c) Let 
$\rho=\alpha_\ivec\,\theta_\ivec\in\Gamm_\ii^{\ot k}$,
$\alpha_\ivec=(\alpha_{i_1}, \dots,\alpha_{i_k})\in\C^k$. Since
$\rho\ot1=\sum\rho_{(0)}\ot\rho_{(1)}$ by  (a)  and (b) we have
\begin{equation} 
\begin{split} S(\rho) 
&=\ts\sum
  S(\rho_{(0)})\rho_{(1)}=S(\alpha_\ivec\,\theta_\jvec)\,
  v^\jvec_\ivec=\alpha_\ivec\, S(\theta_\jvec)v^\jvec_\ivec
\\
&=(-1)^{\half k(k+1)} 
\alpha_\ivec\,\eta_\jinv\, v^\jvec_\ivec =(-1)^{\half k(k+1)} 
\alpha_\ivec\,\sig_\wo(\theta_\ivec)= (-1)^{\half k(k+1)}\sig_\wo(\rho). 
\end{split} 
\end{equation}
(d)  Since $S(\ker\ant)\subseteq \ker\ant$ by \cite[Prop.\,13.10, p.\,489]{KS}
and $\sig_\wo(\ker\ant)\subseteq\ker\ant$ by \rf[bwo], the antipode $S$ and
$\sig_\wo$ are well-defined on the quotient $\gdw[k]$.
\end{proof} 
\begin{lemma} \label{biinv-l} 
Let $\lam\vdash k$, $\lam=\lam'$, be a symmetric partition of $k$ and let
$\pi\in\H$ be a subidempotent of $\pi_\lam$. Then
\begin{equation}    \label{woq} 
\bsig_\wo \pi=(-1)^{\ell(\wo)}\pi. 
\end{equation} 
\end{lemma} 
\begin{proof} (a) We use the Iwahori-Hecke algebra $\hh$ over $\C(q)$. First
we show that   $T_\wo^2z_\lam=z_\lam$, $\lam\vdash k$,
$\lam=\lam^\prime$. By \rf[bwo], $T_{\wo}^2$ is central. Hence there exist
coefficients $\alpha_\lam\in\C(q)$ with $T_\wo^2=\sum_\lam\alpha_\lam
z_\lam$. Since $T_\wo$ is invertible, $T_\wo^{-2}= \sum_\lam\alpha_\lam^{-1}
z_\lam$. By Lemma\,\ref{invzentral}, and $\lam=\lam^\prime$  we obtain,
$\alpha_\lam z_\lam=(\alpha_\lam z_\lam)^\prime=(T_\wo^2
z_\lam)^\prime=z_\lam T_\wo^{-2}=\alpha_\lam^{-1} z_\lam$. Hence,
$\alpha_\lam^2=1$. Let $\alpha_\lam(q)=h(q)^{-1}g(q)$ with polynomials $g$ and
$h$. Since $\alpha_\lam(q)$ takes only two values $1$ and $-1$, at least one
of them, say $c$, appears infinitely often: $\alpha_\lam(q_i)=c$,
$i=1,2,\dots$. Then  the polynomial $f(q)=g(q)-c\cdot h(q)$, has infinitely
many zeros. Hence $f\equiv0$ and  $\alpha_\lam(q)\equiv c$  is
constant. Specialization at $q=1$ yields $T_\wo^2|_{q=1}=1$. Hence, $\alpha_\lam=1$.
\\ 
(b) By definition of $\pi_\lam$, \rf[bwo], $\wo=w_\circ^{-1}$, \rf[relianti],
and (a):
\begin{alignat*}{2} 
\bsig_\wo \pi
    &=\bsig_\wo\pi_\lam \pi= t_\lam^{-1}\bsig_\wo\AAA_k (1\ot z_\lam)\pi\,
                      &=&\,t_\lam^{-1}\AAA_k((-1)^\lo T_{\wo}^{\ast\prime} 
                           \ot T_\wo)(1\ot z_\lam)\pi \nn 
\\ 
&=t_\lam^{-1}(-1)^\lo\AAA_k(1\ot T_\wo^2z_\lam)\pi\,
                         &=&\, t_\lam^{-1}(-1)^\lo\AAA_k(1\ot z_\lam)\pi \nn 
\\ 
&=(-1)^\lo\pi_\lam \pi=
(-1)^\lo \pi.
\end{alignat*} \end{proof} 
\begin{cor} For $\rho\in\gdwi[k]$, 
\begin{equation}        
\label{e-antipode}
\sig_\wo(\rho)=(-1)^{\half k(k-1)}\rho\quad\text{and}\quad S(\rho)=(-1)^k\rho. 
\end{equation} \end{cor} 
\begin{proof} 
Since $\gdwi[k]$ and $\bigoplus_{\lam\in\PP_k}\Mor(1,\phi(\varPi^\tau_\lam))$
are isomorphic linear spaces (by the proof of Theorem~\ref{bi}~(i)),
bi-invariant $k$-forms $\rho$ are in one-to-one correspondence with rank~1
subidempotents $\pi$ of $\sum_{\lam\in\PP_k}\varPi^\tau_\lam$. Since
$\ell(\wo)=\half k(k-1)$ the first part  follows from
Lemma\,\ref{biinv-l}. Combining this  with Proposition\,\ref{antipode} gives $S(\rho)=(-1)^k\rho$.
\end{proof} 
Now we are ready to complete the proof of Theorem\,\ref{bi}\,(ii),  (iii), and
(iv). Using \rf[e-antipode] and \rf[srr]  we obtain for $\rho_1\in\gdwi[k]$
and $\rho_2\in\gdwi[n]$
\begin{alignat*}{2} 
\rho_1\land\rho_2
&=(-1)^{k+n}S(\rho_1\land\rho_2)    &=&(-1)^{k+n+kn}S(\rho_2) \land S(\rho_1)
\\
&=(-1)^{k+n+kn+k+n}\rho_2\land \rho_1&=&(-1)^{kn}\rho_2\land\rho_1. 
\end{alignat*} 
Let $\rho\in\gdwi[k]$. Since the differential $\dd$ commutes with the
antipode, increases the degree of a form  by 1, and maps bi-invariant forms
into bi-invariant forms, again by \rf[e-antipode] we have
\begin{equation}\nn 
\dd \rho=\dd\bigl((-1)^kS(\rho)\bigr)=(-1)^k S(\dd\rho)= 
(-1)^k(-1)^{k+1}\dd\rho=-\dd\rho. 
\end{equation}
Hence, $\dd\rho=0$. Since $\dd$ is linear, each bi-invariant form is closed. 
\\ 
Suppose $\rho_1,\rho_2\in\gdwi[k]$ represent the same de Rham cohomology
class. Then there exists $\rho\in\gdwi[k-1]$ with
$\rho_1-\rho_2=\dd\rho$. Since $\dd$ intertwines the right coaction $\dr$ on
$k-1$ and $k$ forms, the pre-image $W=\dd^{-1}(\dd\rho)$  is  a $\dr$
invariant subspace of $\gdwl[k-1]$. Since  $\dr{\uhr}\langle \dd\rho\rangle$
is the trivial comodule, by Schur's Lemma,  $\dr{\uhr}W$ is a multiple of the
trivial comodule. Consequently, $W\subseteq \gdwi[k-1]$ and
$\dd\rho=0$. Hence, $\rho_1=\rho_2$. The proof of Theorem\,\ref{bi} is complete.
\section{Proof of Theorem\,\ref{usw}}\label{s-univers}
Recall that a {\em quadratic algebra}, see \cite[Sect.\,3]{Ma}, is an $\NO$-graded 
algebra $B=\bigoplus_{i\ge0}B_i$ with $B_0=\C$, $B_1$ 
generates $B$, and the ideal of relations among the elements of $B_1$ is 
generated by a subspace $\Rel(B)\subseteq B_1\ot B_1$. It is convenient to 
write $B=\{B_1,\Rel(B)\}$. The {\em dual}\/ quadratic algebra to 
$B=\{B_1,\Rel(B)\}$ is $B^!=\{B_1^\ast, \Rel(B)^\perp\}$. Here 
$B_1^\ast$ denotes the dual 
vector space to $B_1$ and $\Rel(B)^\perp=\{t\in B_1^\ast\ot B_1^\ast\colon 
t(r)=0, r\in \Rel(B)\}$.
\begin{prp}                              \label{p-dual}
Let $\AA$ be one of the Hopf algebras
$\oglqn$ or $\oslqn$ and let $\Gamm=\Gamm_{-\tz}$. Then $\gdsl$ is
 a quadratic algebra. It is dual to the algebra 
\begin{equation}\nn
B(R^{\tau}):=\C\langle L^i_j\colon i,j=1,\dots,N\rangle/
(L_2\Rda{\tau} L_2\Rda{\tau}-\Rda{\tau} L_2\Rda{\tau}
L_2).
\end{equation}
\end{prp}
\begin{remark}
The Proposition is also valid for bicovariant FODC on quantum groups of types
B, C, and D, when the corresponding $\Rda{}$-matrices are used. The algebra
$B(R)$ is called {\em reflection equation algebra}\/ or {\em algebra of
braided matrices}.
\end{remark}
\begin{proof}
We carry out the proof for
$\Gamm=\Gamm_{+,z}$. The proof for $\Gamm_{-,z}$ is analogous. 
Throughout the proof we sum over repeated indices.
By definition,  $\gdsl=\{\gl, 
\ker(I-\sig_+)\}$ is a  quadratic algebra. Using  a simple argument from linear 
algebra, $(\gdsl)^!=\{\Gamm_\ell^\ast, \im (I-\sig_+)^\fettt\}$.
Let $\{Y_{ij}:i,j=1,\dots,N\}$ be a basis of $\Gamm_\ell^\ast$ dual  to 
the basis $\{\theta^i_j:i,j=1,\dots,N\}$ of $\Gamm_\ell$.
We identify the tensor algebra over $\Gamm_\ell$ and $\Gamm_\ell^\ast$
with the free associative unital complex algebra $\C\langle\theta^i_j\rangle$ and
$\C\langle Y_{ij}\rangle$, respectively. Now we compute the relation subspace
$\im(I-\sig_+^\fettt)$ of $(\gdsl)^!$. By \rf[sig], \rf[e-mor2] and $\Rda[rs]{ab}=\Rda[ab]{rs}$ we obtain
\begin{equation}
\begin{align*}
Y_{ijkl}&:=(I-\sig_+^\fettt)^{prce}_{ijkl}Y_{pr} Y_{ce}
\\
&=Y_{ij}  Y_{kl}
-\Rre[xy]^{jk}(\check{R}^{-1})_{pv}^{ix}\Rda[we]{yl}(\grave{R}^-)^{vw}_{rc}Y_{pr}  Y_{ce}
\\
&=Y_{ij}  Y_{kl}-\Rda[yk]{xj}(\Rdam)_{xi}^{vp}\Rda[we]{yl}(\grave{R}^-)^{vw}_{rc}
Y_{pr}  Y_{ce}
\\
&=Y_{ij}  Y_{kl}-\Rda[xj]{yk}(\Rdam)^{xi}_{vp}Y_{pr}  Y_{ce}\Rda[yl]{we}
q^{-2w+2r}(\Rdam)^{vr}_{wc}.
\end{align*}
\end{equation}
We choose new variables $L^a_b=q^{2b}Y_{ab}$,  multiply the above
equation by $q^{2l+2k}(\Rdam)^{aj}_{bk}$, and sum over $k$ and $j$ (no summation over $a,b,i,l$):
\begin{equation}
\begin{align}
Y_{ijkl}q^{2l+2k}(\Rdam)^{aj}_{bk}
=&\,
q^{-2j+2k}(\Rdam)^{aj}_{bk}L^i_j  L^k_l
-   q^{2l}\Rda[xj]{yk}(\Rdam)^{aj}_{bk}q^{2k-2x}q^{2x}
(\Rdam)_{vp}^{xi}\nn 
\\
&\,\times L^p_z  L^c_e q^{-2e}\Rda[yl]{we}q^{-2w+2r}(\Rdam)^{vr}_{wc}.\nn
\\
\intertext{
Since $\delta_{ax}\delta_{by}=\Rda[xj]{yk}(\Rdam)^{aj}_{bk}q^{2k-2x}$ (no
summation over $x$) 
by \rf[e-mor2] we obtain  in matrix notation $L_2=I\ot L$}
Y_{ijkl}q^{2l+2k}(\Rdam)^{aj}_{bk}
=&\,q^{-2b+2a}L^i_j  (\Rdam)^{aj}_{bk}L^k_l-q^{2l+2a-2e-2w}
(\Rdam)^{ai}_{vp}L^p_r  L^c_e\Rda[bl]{we}(\Rdam)^{vr}_{wc}
\nn \\
=&\, 
q^{-2b+2a}(L_2\Rdam L_2)^{ai}_{bl}-q^{-2b+2a}(\Rdam L_2)^{ai}_{vr}
(\Rdam)^{vr}_{wc}(L_2 \Rda{})^{wc}_{bl}.  \nn
\end{align}
\end{equation}
Multiplying the above equation by $q^{-2a+2b}(\Rdam)^{bl}_{mn}$ and summing over $b,l$ (no summation over $a,i,m,n$) gives
\begin{equation}\nn
L_{aimn}:=Y_{ijkl}q^{2l+2k+2b-2a}(\Rdam)^{aj}_{bk}(\Rdam)^{bl}_{mn}=
(L_2\Rdam L_2\Rdam-\Rdam L_2 \Rdam L_2)^{ai}_{mn}.
\end{equation}
Since $Y_{ijkl}=L_{xiyz}q^{-2l-2v}\Rda[vl]{yz}\Rda[xj]{vk}$ (no summation over $l$), 
$\langle L_{ijkl}:i,j,k,l=1,\dots,N\rangle=\langle
Y_{ijkl}:i,j,k,l=1,\dots,N\rangle$; the  proof is complete.
\end{proof}
It was shown by Majid \cite[Theorem\,7.4.1 and Theorem\,10.3.1]{Maj} 
that the vector space $\AA(R)$  admits {\em another}\/ product 
such that  $\AA(R)$ becomes isomorphic to
the  associative unital algebra $B(R^+)$.
Note that changing the role 
of $\Rda{}$ and $\Rda{-1}$  does not effect the definition of $\AA(R)$.
Hence 
$\AA(R)$, $B(R^+)$,  and $B(R^-)$ are isomorphic  as graded linear spaces.  
It is well known that $\AA(R)$
obeys a linear basis $\{(u^1_1)^{k_{11}}(u^1_2)^{k_{12}} \cdots
(u^N_N)^{k_{NN}}: k_{11},k_{12},\dots,k_{NN}\in\NO\}$, see
\cite[Theorem\,3.5.1]{PW}. Hence $B(R^\tau)$ has a Poincar\' e-Birkhoff-Witt basis
too. By  \cite[Theorem\,5.3]{Priddy}, $B(R^\tau)$ is Koszul. By
\cite[Proposition\,7, Sect.\,9]{Ma} and by Proposition\,\ref{p-dual} we
have $P(\gdsl,t)P(B(R^\tau),-t)=1$. Hence $P(\gdsl,t)=(1+t)^{N^2}$. Note that
$\Js\subseteq\JJ$ by \rf[sym]. Since $P(\gdsl,t)=P(\gdwl,t)$ by
Theorem\,\ref{left}, the ideals $\Js$ and $\JJ$ coincide; the proof of
Theorem\,\ref{usw}\,(i) is complete.
\\
Now we want to compare the universal exterior algebra with the second
antisymmetrizer exterior algebra. It turns out that the bi-invariant 2-form
$\theta^2$ is useful to decide whether or not these two differential Hopf
algebras coincide.
We often need  the following formula. Applying  $\mu(S\ot\id)\dl$, $\mu$
denotes the multiplication,
to equation  \rf[inner] gives
\begin{equation}\label{e-ta}
\theta\rac a=\ve(a)\theta+\om(a).
\end{equation}
The next  lemma is also of interest for its own.
\begin{lemma}\label{l-inner}
Let $\Gamm$ be an inner  bicovariant FODC over $\AA$ with $\dd a=\theta
a-a\theta$, $a\in\AA$. Let $\gd$ be an arbitrary bicovariant differential calculus over
$\AA$ that contains $\Gamm$ as its first order part.
\\
 {\rm (i)} Then we have for $a\in\AA$
\begin{align}
\theta^2\rac a &=\ve(a)\theta^2 +
                 \theta\om(a)+\om(a)\theta-\dd\om(a),
\label{e-tta}
\\
(\dd\theta)\rac a&=\ve(a)\dd\theta +
\theta\om(a)+\om(a)\theta-\dd\om(a).
\label{e-dta}
\end{align}
{\rm (ii)} If $\dd\theta=2\theta^2$ or if $\theta^2=0$, then $\gd$ is inner with $\dd\rho=\theta\rho- (-1)^n\rho\theta $,
$\rho\in\Gamm^n$.  In this situation $\theta^2$ is central in $\gd$.
\end{lemma}
\begin{proof}
By  \rf[e-ta] and \rf[MC] we have
$\theta^2\rac a=
\theta\rac a_{(1)}(\theta\rac
a_{(2)})=(\ve(a_{(1)})\theta+\om(a_{(1)}))(\ve(a_{(2)})\theta+\om(a_{(2)}))=\ve(a)\theta^2+\theta\om(a)+\om(a)\theta-\dd\om(a)$.
Using $\dd\theta\, b=\dd(\theta b)+\theta \dd b$, $\theta b=b\theta +\dd
b$, \rf[e-ta], and  \rf[MC] we obtain
\begin{align*}
\dd\theta\rac a&=\ts\sum Sa_{(1)}\dd \theta a_{(2)}=\ts\sum S a_{(1)}(\dd(\theta
a_{(2)})+\theta \dd a_{(2)})
\\
&=\ts\sum Sa_{(1)}\dd(a_{(2)}\theta+\dd a_{(2)})+\ts\sum Sa_{(1)}\theta a_{(2)}Sa_{(3)}\dd a_{(4)}
\\
&=\ts\sum S a_{(1)}(\dd a_{(2)}\theta +a_{(2)}\dd\theta)+
\ts\sum(\ve(a_{(1)})\theta+\om(a_{(1)}))\om(a_{(2)})
\\
&=\om(a)\theta+\ve(a)\dd\theta+\theta
\om(a)-\dd\om(a).
\end{align*}
(ii) Suppose first $\dd\theta=2\theta^2$. Taking the difference
$2\times\rf[e-tta]-\rf[e-dta]$, we have $0=\theta\om(a)+\om(a)\theta
-\dd\om(a)$. Since each $\rho\in\Gamm^1$ can be written in the form
$\rho=\sum_ib_i\om(a_i)$, $a_i,b_i\in\AA$, $\dd\rho=\sum_i((\theta
b_i-b_i\theta)\om(a_i)+b_i(\theta\om(a_i)+\om(a_i)\theta)=\theta\rho+\rho\theta$.
The proof for $\rho\in\Gamm^k$, $k\ge2$, goes by induction on $k$.
We show $\theta^2$ is central.
Inserting $\dd\om(a)=\theta\om(a)+\om(a)\theta$ into \rf[e-tta] gives
$\theta^2\rac a=\ve(a)\theta^2$. Equivalently, $\theta^2 a=a\theta^2$. Since
$\Gamm$ is inner, $\theta^2\rho=\rho\theta^2$, $\rho\in\Gamm^k$, follows by
induction on $k$. Suppose now $\theta^2=0$. Let $\rho=\sum_i a_i\dd b_i$ be a $1$-form. Then $\dd\rho=\sum_i\dd a_i\dd
b_i= \sum_i(\theta a_i-a_i\theta)\dd b_i=\theta\rho-\sum_i a_i\theta(\theta
b_i-b_i\theta)=\theta\rho+\sum_i a_i(\theta b_i-b_i\theta)\theta+\sum_i
a_i(-\theta^2 b_i+ b_i\theta^2)=\theta\rho+\rho\theta$. The proof for
$\rho\in\Gamm^k$, $k\ge2$, goes by induction on $k$.
\end{proof}
Throughout the remainder let
$\eta=\SS(U)=\sum_{i,j}q^{-2i}\om_{ij}\ott\om_{ji}$. Note that
$(\Gamm^{\ot 2})_\ii=\langle\vt\ot\vt,\eta\rangle$ since
$\dim\Mor(1,u^\fettc\ott u\ott u^\fettc\ott u)=2$.
Using  \rf[e-trans] and \rf[e-ta], similarly to the  calculations
in the preceding lemma one gets
\begin{align}    
\vt\rac a&=\ntz\om(a)+\ve(a)\vt, \label{e-vta}
\\
(\vt\ot\vt)\rac b&=\ntz^2\SS(b)+\ntz\bigl(\om(b)\ot\vt+\vt\ot\om(b)\bigr)+\ve(b)\vt\ot\vt.
\label{e-vtvtb}
\\
\intertext{By  the definition of $\SS$, $\om(\wt{a}b)=\om(a)\rac b$,
\rf[e-vta], $\om(U)=\vt$, and $\ve(U)=\fr{s}$ one easily checks}        
\SS(ab)&=\SS(a)\rac b+(\om(a)\rac b_{(1)})\ot \om(b_{(2)})+\om(b_{(1)})\ot 
(\om(a)\rac b_{(2)})+\ve(a)\SS(b),
 \label{e-sab}
\\
\SS(Ub)&=\eta\rac b+(2\ntz
+\fr{s})\SS(b)+\vt\ot\om(b)+\om(b)\ot\vt.\label{e-sub}
\end{align}
\begin{prp}    \label{su}
Let $\AA$ be one of the Hopf algebras $\oglqn$ or $\oslqn$,
$\Gamm=\Gamm_{\tau,z}$, and $q$ be transcendental.  
Then $\gdu\cong\gds$ as differential Hopf algebras if and only if
$\theta^2=0$ in $\gdu$.
\end{prp}
\begin{proof}
 $\rightarrow$ Since $\theta$ is bi-invariant, $(I-\sig)(\theta\ota\theta)=0$.
Hence $\theta^2=0$ in $\gds$.
\\
$\leftarrow$ Suppose  $\theta^2=0$ in $\gdu$.
By universality of $\gdu$, $\gds$ is a quotient of $\gdu$. 
To the converse relation.  
In \cite[Sect.\,5.3]{LS} it was shown that the quadratic algebra $\gdsl$ has 
defining relations
\begin{equation} \label{plus}
\Rda{\tau}\vN_2\Rda{\tau}\vN_2\Rda{\tau}
+\vN_2\Rda{\tau}\vN_2  =0.
\end{equation}
We have to prove that $\vN=(\theta^i_j)$, $\theta^i_j\in\gdu$, satisfies \rf[plus]. Then $\gdu$ would be a  quotient of $\gds$,
and both differential Hopf algebras coincide.
\\ 
Inserting $a=u^i_j$ into \rf[e-tta] and using
$-\dd\om(u^i_j)=\sum_x\om_{ix}\om_{xj}$ (by \rf[MC]), we have
$\theta\om_{ij}+\om_{ij}\theta+\sum_x\om_{ix}\om_{xj}=0$, $i,j=1,\dots,N$. In
matrix notation $\vO=(\om_{ij})$, $\theta\vO+\vO\theta+\vO\vO=0$.
By \rf[e-trans], $\theta^2=0$, and
$z^2r_\tau\ne0$  we  obtain
\begin{equation}\label{x}
r_\tau\vN\vN+\theta\vN+\vN\theta=0.
\end{equation}
We complete the proof for the sample $\Gamm=\Gamm_{+,z}$. 
The proof for $\Gamm=\Gamm_{-,z}$ is analogous and uses $\Rda{-2}=-\QM
\Rda{-1}+I$ instead of $\Rda{2}=\QM\Rda{}+I$.
Applying $\,\,\rac u$ to \rf[x], using  $\Rda{2}=\QM\Rda{}+I$,
and again \rf[x] (to the underlined terms) we have
\begin{align}         
0=&(r_+\vN_{1}\vN_{1}+\theta\vN_{1}+\vN_{1}\theta)\rac u_{2}=
(r_+ \vN_{1}+\theta I_{1})\rac u_{2}(\vN_{1}\rac u_{2}) +(\vN_{1}\rac u_{2})
(\theta I_{1}\rac u_{2})
 \nn\\
=&(r_+z\Rda{}\vN_{2}\Rda{}+r_+z\vN_{2}+z\theta)z\Rda{}\vN_{2}\Rda{}
+r_+z^2\Rda{}\vN_{2}\Rda{}\vN_{2}+z^2\Rda{}\vN_{2}\theta\Rda{}
\nn\\
=& r_+z^2\Rda{}\vN_{2}(\QM\Rda{} +\ul{I})\vN_{2}\Rda{}+r_+z^2\vN_{2}\Rda{}\vN_{2}\Rda{}+
\ul{z^2\Rda{}(\theta\vN_{2}+\vN_{2}\theta)\Rda{}}+r_+z^2\Rda{}\vN_{2}\Rda{}\vN_{2}
\nn\\
=&r_+z^2(\QM\Rda{}\vN_{2}\Rda{}\vN_{2}\Rda{}+\vN_{2}\Rda{}\vN_{2}\Rda{}+
\Rda{}\vN_{2}\Rda{}\vN_{2})=
r_+z^2\Rda{} (\Rda{}\vN_{2}\Rda{}\vN_{2}\Rda{}+\vN_{2}\Rda{}\vN_{2}).\nn
\end{align}
Since $r_+z^2\Rda{}$ is invertible, $\vN$ satisfies 
\rf[plus]; this completes the proof.
\end{proof}
Before proving $\theta^2=0$ for all cases except for $N=2$, $\tau=+$, and 
$z^2=\qm[2]$, we need a rather technical lemma that   describes the right ideal
$\RR_{\tau,z}$, $\ntz\ne0$.  First  recall the defining parameters (for
$\ntz,r_\tau$ and $\fr[\pm]{r}$ see before \rf[lin]; for  $\fr[\pm]{s}$ see
before  \rf[e-mor1]). We use the same notations as in \cite{SS2,SS3}. For $\eps\in\{+,-\}$ define
$(v_\eps)^i_j:=\sum_{k,n,m}q^{-2k}(P_\eps)^{ki}_{nm}u^n_ku^m_j,$ and
$V_\eps:=\sum_{i}q^{-2i}(v_\eps)^i_i$.
 Set  $\nu_\tz=\ntz^{-1}(q^{-2\tau}z^{-N}-1)$ and 
\begin{xalignat*}{3} 
\fr[\tau]{s}^+&=\fr{s}+q^\tau   r_\tau\QP,
&\fr[\tau]{s}^-&=\fr{s}+q^{-\tau}r_\tau\QP,&
\\
 \fr[+]{t}&=\fr{s}+1, &\fr[-]{t}&=\fr{s}-q^{-2N},&
\\
\alpha_\eps&=\QP^2\fr[\eps]{s}^{-1}, &\gamma_\eps&=q\QP\fr{s}^{-1}
\fr[\eps]{t}^{-1},&
\\
\lam^+_{\tau,z}&=\QP^{-1}q^{\tau} z\fr[+]{s},
&
\lam^-_{\tau,z}&=\QP^{-1}q^{-\tau} z\fr[-]{s},&
\\
   \mu^+_{\tau,z}&= \fr[+]{t}(z^2\fr[\tau]{s}^+-\fr{s})(q\QP \ntz)^{-1},
  & \mu^-_{\tau,z}&=\fr[-]{t}(z^2\fr[\tau]{s}^--\fr{s})(q\QP \ntz)^{-1}.&
\end{xalignat*}
The following identities are easily checked
\begin{align}
q^{2\tau}\fr[+]{s}+q^{-2\tau}\fr[-]{s}&=\fr[\tau]{s}^++\fr[\tau]{s}^-=2\fr{s}+r_\tau\QP^2,
\label{est1}\\
q^{\tau}\fr[+]{s}+q^{-\tau}\fr[-]{s}&=\QP(\fr{s}+r_\tau).
\label{est2}
\end{align}
By \rf[est2], $\lam_{\tau,z}:=\lam^+_{\tau,z}+\lam^-_{\tau,z}=z\fr[\tau]{r}$.
Define the $N^2{\times}N^2$-matrices $P_\iota=(P_\iota)^{ab}_{rs}$,
$\iota=0,1$, by $(P_0)^{ab}_{rs}=\frac{1}{\fr{s}}
q^{-2a}\delta_{ab}\delta_{rs}$, $P_1=I-P_0$. Obviously, $P_0$ and $P_1$ are
projection operators and span $\Mor(u^\fettc\ott u)$. 
For $(\eps,\iota)\in\{+,-\}\times\{0,1\}$ define the map
$Q_{\eps,\iota}\in\End(V^{\ot 2},V^{\ot 4})$ by 
\begin{equation*}
(Q_{\eps,\iota})^{abmn}_{rs}
=\sum_{x,y,z}(P_\eps)^{ab}_{xy}(P_{\iota})^{mn}_{yz}(P_{\eps})^{xz}_{rs}.
\end{equation*}
Since  $P_\eps\in\Mor(u\ott u)$ and 
$P_\iota\in\Mor(u^\fettc\ott u)$, 
$Q\in\Mor(u\ott u,u\ott u\ott u^\fettc\ott u)$ is easily checked.
\begin{lemma}\label{l-ql} 
Let $\Gamm=\Gamm_{\tau,z}$. 
{\rm (i)} 
\begin{equation}                    \label{e-uu}
u^a_ru^b_s-\delta_{ar}\delta_{bs}\equiv \sum_{\eps,\iota,m,n}
c_{\eps,\iota}(Q_{\eps,\iota})^{abmn}_{rs}
\wt{u^m_n}\mod\RR_{\tau,z},
\end{equation}
where $c_{\eps,1}=\alpha_\eps\lam^\eps_{\tau,z}$ and
$c_{\eps,0}=\fr{s}\gamma_\eps\mu^\eps_{\tau,z}$.
\\
{\rm (ii)}
\begin{equation*}
\SS(u^a_ru^b_s)=\sum_{\eps,\iota,\iota'}c_{\eps,\iota}c_{\eps,\iota'}
(Q_{\eps,\iota})^{abmn}_{xy}(Q_{\eps,\iota'})^{xyvw}_{rs}\om_{mn}\ot\om_{vw},
\end{equation*}
(summation over $m,n,x,y,v$, and $w$).
\\
{\rm (iii)} $\SS(W_\eps)=
A^\eps_\tz\eta
+B^\eps_\tz\vt\ot\vt$, where $W_\eps=V_\eps-\mu^\eps_\tz U$ and 
\begin{equation} \label{abpm}
A^\eps_\tz=\alpha_\eps(\lam^\eps_{\tau,z})^2-\mu^\eps_\tz,\quad
B^\eps_\tz=\gamma_\eps(\mu^\eps_{\tau,z})^2-\fr{s}^{-1}
\alpha_\eps(\lam_\tz^\eps)^2.
\end{equation}
\end{lemma}
\begin{proof}
(i)  Throughout the proof 
we sum over repeated indices.  Recall 
\cite[formula (2),\,(3), p.\,318]{SS3}:
$(v_\pm)^i_j-\frac{1}{\fr{s}}\delta_{ij}V_\pm\equiv\lam^\pm_{\tau,z}(u^i_j-\frac{1}{\fr{s}}\delta_{ij}U)\mod\RR_{\tau,z}$
and 
\begin{equation} \label{e-vpm}
\wt{W_\eps}=\wt{V_\eps}-\mu^\eps_{\tau,z}\wt{U}
\equiv0\mod\RR_{\tau,z}.
\end{equation}
Inserting this into the main
formula of the proof of \cite[Lemma\,4.4]{SS2}
\begin{equation*}
\begin{split}
u^a_ru^b_s-\delta_{ar}\delta_{bs}&\equiv
\fr[+]{s}^{-1}\QP^2(P_+)^{ab}_{xy}
\wt{v_+}^y_z(P_+)^{xz}_{rs}+\fr[-]{s}^{-1}\QP^2(P_-)^{ab}_{xy}
\wt{v_-}^y_z(P_-)^{xz}_{rs}
\\
&\,\,
-q\QP\fr[+]{s}^{-1}\fr[+]{t}^{-1}(P_+)^{ab}_{rs}\wt{V_+}
-q\QP\fr[-]{s}^{-1}\fr[-]{t}^{-1}(P_-)^{ab}_{rs}\wt{V_-},
\end{split}
\end{equation*} 
and using 
$\fr[\eps]{s}=q^{-1}\QP\,\fr[\eps]{t}-\fr{s}$, a number of long computations
gives \rf[e-uu].
\\
(ii) By 
the definition of $\SS$, \rf[e-uu], and  $\om(a)=\om(b)$ for $a\equiv b\mod 
\RR_{\tau,z}$, 
\begin{equation}\nn
\SS(u^a_ru^b_s)=\om(u^a_xu^b_y)\ot\om(u^x_ru^y_s)=
\sum_{\eps,\eps'\iota,\iota'}
c_{\eps,\iota}c_{\eps',\iota'}(Q_{\eps,\iota})^{abmn}_{xy}(Q_{\eps',\iota'})
^{xyvw}_{rs} \om_{mn}\ot\om_{vw}.
\end{equation}
Since $P_\eps P_{\eps'}=\delta_{\eps,\eps'}P_\eps$, (ii) follows as claimed.
\\
(iii) By definition of $V_\eps$, $P_\eps P_{\eps'}=\delta_{\eps,\eps'}P_\eps$,
and $q^{-2x-2y}(P_\eps)^{ij}_{xy}=q^{-2i-2j}(P_\eps)^{ij}_{xy}$ (no summation),
\begin{equation}\nn
\begin{split}
\SS(V_\eps)
&=q^{-2i-2j}(P_\eps)^{ij}_{xy}\sum_{\eps',\iota,\iota'}c_{\eps'\iota}c_{\eps'\iota'}
(P_{\eps'})^{xy}_{ab}(P_\iota)^{mn}_{bc}(P_{\eps'})^{ac}_{ed}(P_{\iota'})^{vw}_{df}
(P_{\eps'})^{ef}_{ij}\om_{mn}\ot\om_{vw}
\\
&=q^{-2x-2y}\sum_{\iota,\iota'}c_{\eps\iota}c_{\eps\iota'}(P_\eps)^{xy}_{ab}(P_\iota)^{mn}_{bc}(P_\eps)^{ac}_{xd}(P_{\iota'})^{vw}_{dy}\om_{mn}\ot\om_{vw}.
\end{split}
\end{equation}
Note that $T_\eps=(T_\eps^{bcdy})$,
$T_\eps^{bcdy}=q^{-2x-2y}(P_\eps)^{xy}_{ab}(P_\eps)^{ac}_{xd}$ (no summation 
over $y$), belongs to
$\Mor(1,u^\fettc\ott u\ott u^\fettc\ott u)$. Since there exists a linear
isomorphism of the latter space onto
$\Mor(u^\fettc\ott u)$, it is 
two-dimensional.  The two  morphisms $Q_0$ and $Q_1$, $Q_\iota
^{bcdy}:=q^{-2y}(P_\iota)^{bc}_{yd}$, $\iota=0,1$ (no summation over $y$), form a basis 
of  $\Mor(1,u^\fettc\ott u\ott u^\fettc\ott u)$.
By \rf[e-mor2], $T_\eps=\alpha_\eps^{-1}Q_1+\fr{s}^{-1}\gamma_\eps^{-1}Q_0$.
Moreover,  $(P_\iota\ot P_{\iota'})
Q_{\iota''}= \delta_{\iota,\iota'}\delta_{\iota,\iota''}Q_\iota$,
 $Q_0^{mnvw}\om_{mn}\ot\om_{vw}=\frac{1}{\fr{s}}\vt\ot\vt$, 
and  $Q_1^{mnvw}\om_{mn}\ot\om_{vw}=
 \eta-\frac{1}{\fr{s}}\vt\ot\vt$.
Hence,
\begin{equation}
\begin{split}
\SS(V_\eps)&=(\alpha_\eps^{-1}c_{\eps,1}^2
Q_1^{mnvw}+\fr{s}^{-1}\gamma_\eps^{-1}c_{\eps,0}^2Q_0^{mnvw})
\om_{mn}\ot\om_{vw}
\\
&=\alpha_\eps(\lam_\tz^\eps)^2\eta+\bigl(\gamma_\eps
(\mu_\tz^\eps)^2-\fr{s}^{-1}\alpha_
\eps(\lam_\tz^\eps)^2\bigr)\vt\ot\vt.
\end{split}
\end{equation}
Since $\SS(1)=0$ and $\SS(U)=\eta$, the proof of (iii) is complete.
\end{proof}
\begin{lemma}\label{l-t2}
Let $\AA$ be one of the Hopf algebras $\oglqn$ or $\oslqn$. Let $\Gamm=\gtz$, 
$\RR=\RR_\tz$, and suppose $\ntz\ne0$. Then for
$N\ge3$, $\tau\in\{+,-\}$, 
and for $N=2$, $\tau=+$, and $z^2\ne \qm[2]$  we have $\vt^2=0$ in $\gdu$.
\end{lemma}
\begin{proof}
The proof divides into two parts. The easy part is $z^{2N}\ne
q^{-4\tau}$. Treating the critical values $z^{N}=q^{-2\tau}$ and
$z^{N}=-q^{-2\tau}$  is more difficult. Note  that the condition
$z^{N}=q^{-2\tau}$ is necessary (!) for $\AA=\oslqn$, see before \rf[gpmz]. 
To make  notations simpler, throughout the proof we skip the index $\tz$ in
$A^\eps_\tz, B^\eps_\tz,\ntz,\nu_\tz,
\lam^\eps_\tz,\mu^\eps_\tz$, and $\lam_\tz$. 
\\
{\em Case 1}. $ q^{4\tau}z^{2N}\ne1$. 
By \cite[formula (3) p.\,318]{SS3}, and by  \cite[Theorem\,2.1\,(iii)]{SS3}
\begin{equation}\nn
\wt{\dett}\equiv \nu \wt{U}\mod\RR \quad\text{and}\quad
 \wt{\DD}\equiv
-z^Nq^{2\tau}\nu \wt{U}\mod \RR.
\end{equation}
Since $\DD$
and $\dett$ are group-like, by $\SS(\wt{U})=\eta$, and by  
$\om(\wt{U})=\vt$,
\begin{equation}\label{e-xu}
\begin{split}
\SS(\wt{\dett}-\nu \wt{U})
&=\nu ^2\vt\ot\vt-\nu \eta, 
\\
\SS(\wt{\DD}+z^Nq^{2\tau}\nu \wt{U})&=
z^{2N}q^{4\tau}\nu ^2\vt\ot\vt+
z^Nq^{2\tau}\nu \eta.
\end{split}
\end{equation}
Since $ q^{2\tau}z^{N}\ne1$, $\nu \ne 0$; hence 
$\eta=\nu \vt^2=-z^{N}q^{2\tau}\nu \vt^2$ in
$\gdu$. Thus  $\vt^2=0$ since $q^{2\tau}z^N\ne
-1$.
\\
{\em Case 2}. $(z^Nq^{2\tau}-1)(z^Nq^{2\tau}+1)=0$. 
\\
We set $x:=\wt{U}U-\lam \wt{U}$. Note that
$x=\wt{W}_++\wt{W}_-$. By \rf[e-vpm] and by $\mu^+ +\mu^- =\lam +\fr{s}$, 
$x\equiv0\mod \RR $. By Lemma\,\ref{l-ql}\,(iii), $\SS(x)=A\eta+B\vt\ot\vt$ where
$A=A^+ +A^- $ and  $B=B^+ +B^- $.
\\
Inserting $b=U$ into \rf[e-vtvtb] and $b=\wt{U}$ into \rf[e-sub]  the preceding gives
\begin{equation}\label{action}
(\vt\ot\vt)\rac U=\fr{n}^2\eta+(2\fr{n}+\fr{s})\vt\ot\vt,
\quad \eta\rac U=(A -\fr{n}+\fr{s})\eta +(B -2)\vt\ot\vt.
\end{equation}
By \rf[e-sab], $\SS(rb)=\SS(r)\rac b$, $r\in\RR $, $b\in \AA$. In
particular, for $r=x$ and $b=U$ it follows from
\rf[action] that 
\begin{align}
\SS(xU)&=A \eta\rac U+B (\vt\ot\vt)\rac U
\nn\\
&=(A (A -\fr{n}+\fr{s})+B  \fr{n}^2)\eta
+(A (B -2)+B (2\fr{n}+\fr{s}))\vt\ot\vt
\label{xxx}\\
&=:A ^\prime\eta +B ^\prime\vt\ot\vt.\nn
\end{align} 
Suppose the matrix
$\begin{pmatrix}A &B \\A ^\prime&B ^\prime\end{pmatrix}$ is
regular. Then $\vt\ott\vt$ is a linear combination of $\SS(x)$ and
$\SS(xU)$. Hence $\vt\ot\vt\in\SS(\RR )$ as desired. Otherwise the
matrix is singular,  equivalently
$(A -\fr{n} B )(2A -\fr{n} B )=0$. Repeating the above argumentations
with $\{\SS(x),\SS(\wt{W}_+U)\}$, $\{\SS(x),\SS(\wt{W}_-U)\}$, 
$\{\SS(\wt{W}_+),
\SS(\wt{W}_+U)\}$, and $\{\SS(\wt{W}_-),
\SS(\wt{W}_-U)\}$ instead of $\{\SS(x),\SS(xU)\}$, we end up with
$\vt\ott\vt\in\SS(\RR )$ or two possibilities: $2A^+=\fr{n} B^+, 2A^-=\fr{n} B^-$ or
$A^+=\fr{n} B^+, A^-=\fr{n} B^-$.
\\
{\em Case~2.1}. $2A^\eps =\fr{n} B^\eps$, $\eps\in\{+,-\}$.
Suppose first $N\ge3$. 
Since $q$ is transcendental $\lam^\eps $ is nonzero. 
 Inserting the values \rf[abpm] for $A ^\eps$ and $B ^\eps$,
using $\gamma_\eps\mu ^\eps=\fr{s}^{-1}\fr{n}^{-1}
(z^2\fr[\tau]{s}^\eps-\fr{s})$ and $\lam =\fr{s}+\fr{n}$ we
have
\begin{align}
2\alpha_\eps(\lam ^\eps)^2-2\mu^\eps
&=\fr{n}\frac{z^2\fr[\tau]{s}^\eps-\fr{s}} {\fr{s}\fr{n}}
\mu^\eps-\frac{\lam -\fr{s}} {\fr{s}} \alpha_\eps(\lam ^\eps)^2,
\nn\\
(\lam +\fr{s})\alpha_\eps(\lam ^\eps)^2&=(z^2\fr[\tau]{s}^\eps+\fr{s})\mu^\eps.
\label{2.1}
\end{align}
Suppose $\mu^\eps =0$ for some $\eps\in\{+,-\}$. By \rf[2.1],
$\lam +\fr{s}=z\fr[\tau]{r}+\fr{s}=0$, because $\alpha_\eps(\lam ^\eps)^2\ne0$.
By the definition of $\mu^\eps $,
$z^2\fr[\tau]{s}^\eps-\fr{s}=0$. Thus $z=-q^{-\eps\tau}$ and  $q$ is a root
 of unity. This contradicts our assumption. Hence $\mu^+ \mu^- \ne0$.
Inserting $\lam +\fr{s}=\mu ^++\mu ^-$ into \rf[2.1] yields
\begin{align*}
\mu^+ (\alpha_+(\lam^+ )^2
-z^2\fr[\tau]{s}^+-\fr{s})&=-\alpha_+(\lam^+ )^2\mu^{-}, 
\\
\mu^- (\alpha_-(\lam^- )^2
-z^2\fr[\tau]{s}^--\fr{s})&=-\alpha_-(\lam^- )^2\mu^{+}.
\end{align*}
We take the product of  both equations, divide by $\mu^+\mu^-$ and afterwards insert
$\alpha_\eps(\lam ^\eps)^2=q^{2\eps\tau}z^2\fr[\eps]{s}$. This gives
\begin{align*}
\alpha_+(\lam ^+)^2(z^2\fr[\tau]{s}^-+\fr{s})+\alpha_-(\lam ^-)^2
(z^2\fr[\tau]{s}^++\fr{s})&=(z^2\fr[\tau]{s}^-+\fr{s})
                           (z^2\fr[\tau]{s}^-+\fr{s}),
\\
z^2\bigl(q^{2\tau}\fr[+]{s}(z^2\fr[\tau]{s}^-+\fr{s})+
q^{-2\tau}\fr[-]{s}(z^2\fr[\tau]{s}^++\fr{s})\bigr)
&=
z^4\fr[\tau]{s}^+\fr[\tau]{s}^-
-z^2\fr{s}(\fr[\tau]{s}^++\fr[\tau]{s}^-)+\fr{s}^2.
\end{align*}
By \rf[est1] the coefficient of $z^2$ vanishes. Inserting
$\fr[\tau]{s}^\eps=\fr{s}+q^{\tau\eps}r_\tau\QP$ yields
\begin{align*}
0&=z^4(q^{2\tau}\fr[+]{s}(\fr{s}+q^{-\tau}r_\tau\QP)
+q^{-2\tau}\fr[-]{s}(\fr{s}+q^\tau r_\tau\QP)-(\fr{s}+q^{-\tau}r_\tau\QP)
(\fr{s}+q^\tau r_\tau\QP))-\fr{s}^2,
\\
0&=z^4((q^{2\tau}\fr[+]{s}+q^{-2\tau}\fr[-]{s})\fr{s}+(q^{\tau}\fr[+]{s}+q^{-\tau}\fr[-]{s})\QP
r_\tau-\fr{s}^2-\fr{s} r_\tau\QP^2-r_\tau^2\QP^2).
\end{align*}
Using \rf[est1], \rf[est2], and $\fr{s}\ne0$ we have $z^4(\fr{s}+\QP^2 r_\tau)-\fr{s}=0$.
Since $z^{2N}=q^{-4\tau}$, this gives  for $\tau=+$ 
\begin{align*}
(q^{2N+2}+q^{2N}+q^{2N-6}+q^{2N-8}+\cdots+1)^N-q^8(q^{2N-2}+q^{2N-4}
+\cdots+1)^N&=0,
\\
\intertext{and for $\tau=-$}
(q^{2N+2}+q^{2N}+\cdots+q^8+ q^2+1)^N-q^{4N-8}(q^{2N-2}+q^{2N-4}
+\cdots +1)^N&=0.
\end{align*}
This contradicts our assumption that $q$ is transcendental. Hence
$2A^\eps =\fr{n} B^\eps $ is impossible for $N\ge3$. 
Now let $N=2$, $\Gamm=\Gamm_{+,z}$,
$z^2=-q^{-2}$.  Corresponding to
$z=\pm\ii \qm$ ($\ii=\sqrt{-1}$),  the  defining parameter are
$\lam =\lam^+ =\pm\ii(\qm+\qm[5])$, $\alpha_+\lam =\pm\ii(q+\qm)$,
$\mu ^+=-\fr{n}^{-1}(1+\qm[2]+
\qm[4])(\qm[2]+\qm[6])$,
$\fr{n}=\pm\ii\qm-\qm[2]-\qm[4]\pm\ii\qm[5]$,
$\gamma_+(\mu ^+)^2=\fr{n}^{-2}(q^4+q^2+1)(\qm[2]+\qm[6])^2$. Inserting these
values  into \rf[abpm], elementary transformations of equation $2A ^+-\fr{n}B
^+=0$ lead in  both cases to  $(q^2+1)^2(q^4+1)=0$. This contradicts our
assumption that $q$ is transcendental.
\\
{\em Case~2.2}. $A^\eps =\fr{n} B^\eps $, $\eps\in\{+,-\}$.
Similarly to the derivation of \rf[2.1] we obtain $\lam_{}\alpha_\eps
(\lam_{}^\eps)^2=z^2\mu^\eps_{}\fr[\tau]{s}^\eps$. Inserting $\alpha_\eps(\lam^\eps_{})^2=q^{2\eps\tau}z^2\fr[\eps]{s}$ and
dividing by $z^2$ we get
$q^{2\eps\tau}\lam_{}\fr[\eps]{s}=\mu^\eps_{}\fr[\tau]{s}^\eps$. Using $q^{2\tau\eps}\fr[\eps]{s}=\tau\eps\QM^{-1}\QP
r_\tau+\fr[\tau]{s}^\eps$, subtracting $\lam_{}\fr[\tau]{s}^\eps$, and
multiplying by  $\fr[\tau]{s}^{-\eps}$ gives
\begin{align*}
\lam_{}\tau\QM^{-1}\QP r_\tau\fr[\tau]{s}^-&=
\fr[\tau]{s}^+\fr[\tau]{s}^-(\mu_{}^+-\lam_{}),
\\
-\lam_{}\tau\QM^{-1}\QP r_\tau\fr[\tau]{s}^+&=
\fr[\tau]{s}^+\fr[\tau]{s}^-(\mu_{}^--\lam_{}).
\end{align*}
Taking the sum of both equations, using $\mu^+_{}+\mu^-_{}=\lam_{}+\fr{s}$,
and $\fr[\tau]{s}^--\fr[\tau]{s}^+=-p\tau\QM\QP r_\tau$  we obtain
$-\lam_{} \QP^2 r_\tau^2=(\fr{s}-\lam_{})\fr[\tau]{s}^+\fr[\tau]{s}^-$.
Hence $\lam_{}=\ds\frac{\fr{s} \fr[\tau]{s}^+\fr[\tau]{s}^-}
{(\fr[\tau]{s}^+\fr[\tau]{s}^--\QP^2 r_\tau^2)}$. Using
$\fr[\tau]{s}^+\fr[\tau]{s}^--\QP^2 r_\tau^2=\fr{s}(\fr{s}+ r_\tau \QP^2)$ gives
$\lam_{}=\ds\frac{\fr[\tau]{s}^+\fr[\tau]{s}^-}{\fr{s}+ r_\tau\QP^2}$.
 Inserting
$\lam_{}=z(\fr{s}+r_\tau)$  we 
finally have
\begin{equation*}
z(\fr{s}+r_\tau)(\fr{s}+r_\tau\QP^2)=(\fr{s}+q^\tau\QP r_\tau)(\fr{s}+
q^{-\tau}\QP r_\tau).
\end{equation*}
Since $z^{2N}=q^{-4\tau}$, we conclude that $q$ is algebraic.
This contradicts our assumption; the case $A^\eps =\fr{n} B^\eps $ is
impossible.
Summarising the results of Case~2.1 and Case~2.2, $\vt^2=0$ in $\gdu$,  and the proof is 
complete.
\end{proof}
Combining the preceding with Proposition\,\ref{su} and using
$\theta^2=\fr{n}^{-2}\vt^2$ proves 
Theorem\,\ref{usw}\,(ii). 
\\
To prove  Theorem\,\ref{usw}\,(iii) we will show that
$(\ker\ant[2])_\ll=\SS(\RR)\oplus \langle\vt\ot\vt\rangle$. Since
$(\Gamm^\ot)_\ll$ is a free associative algebra, this proves $\vt\ot\vt\not\in\Ju$.
In Step~1 we will demonstrate $\vt\ott\vt{\not\in}\SS(\RR)$. 
In Step~2 we will show that 
$\dim\SS(\RR)\ge9$. Since $\dim(\ker\ant[2])_\ll=10$, by Theorem\,\ref{left},
and since $\SS(\RR)\subseteq(\ker\ant[2])_\ll$ this proves the
assertion. The remainder of this section is exclusively concerned with
the calculus $\Gamm_{+,z}$, $z^2=q^{-2}$, on $\oglqn[2]$ or $\oslqn[2]$, so we
skip the index ${+,z}$.
\\
{\em Step~1}. We fix $\AA=\oglqn[2]$ and discuss the 
the necessary modifications for
$\AA=\oslqn[2]$ parallel. Roughly speaking, all conclusions remain valid 
if we formally set $\DD=\dett=1$.
Let us collect the defining 
parameters. Since $\nu=0$, $\wt{\dett}$
and $\wt{\DD}$ belong to $\RR$, and furthermore, by \rf[e-xu], 
$\SS(\dett)=\SS(\DD)=0$.
We have $\lam^-=\mu^-=0$. 
Corresponding to $z=\pm \qm$ we find $\fr{n}=\pm \qm-\qm[2]-\qm[4]\pm\qm[5]$ and
\begin{align*}
\lam^+&=\pm(\qm+\qm[5]),&c_{+,1}&=\pm\QP,
\\
\mu^+&=\pm\qm[5](q\pm1)^2(q^2\mp q+1),&c_{+,0}&=\pm \QP
(q\pm 1)^2(q^2\pm q+1)^{-1},
\\
A&=\qm[6](q\mp1)(q^5\mp1),&B&=\pm2\qm(q^5\mp1)(q^3\mp1)^{-1}.
\end{align*}
Our aim is to show that $\xi=\SS(x)$, $x=\wt{U}U-\lam\wt{U}$, is the only bi-invariant element of
$\SS(\RR)$.
Since $\AA$ is cosemisimple and since $\SS\colon\RR\to\gl\ot\gl$
 intertwines 
$\adr$ and  $\dr$, $\ker(\SS{\uhr}\RR)$ has
an $\adr$-invariant complement $\RR_1$ in $\RR$ where $\SS$ acts
injectively. Thus if $\SS(r)$, $r\in\RR$, is bi-invariant,
$r_1=(\SS{\uhr}\RR_1)^{-1}(\SS(r))$ satisfies  $\SS(r_1)=\SS(r)$ and $\adr(r_1)=r_1\ot
1$. Consequently,
the subspace of bi-invariant elements of $\SS(\RR)$ equals 
$\SS(\RR_\inv)$, where $\RR_\inv=\RR\cap\AA_\inv$ and
$\AA_\inv=\{a\in\AA:\adr a=a\ot1\}$.  
\\
We are going to describe $\AA_\inv$. Let $v=(v^i_j)$ be an irreducible
corepresentation of $\AA$. Since $\adr v^i_j=\sum v^x_y\ot (v^\fettc\ott
v)^{xy}_{ij}$, by Schur's Lemma there is a unique up to scalars 
$\adr$-invariant element $\zeta(v)$, $\zeta(v)\in\langle v^i_j\rangle$.
On the other hand since  $\AA$ is cosemisimple,
$\AA_\inv=\bigoplus_v\langle\zeta(v)\rangle$, where we sum over all irreducible
corepresentations $v$ of $\AA$. The set of irreducible corepresentations of
$\oglqn[2]$ and  $\oslqn[2]$  has been  described in 
\cite[Sects.\,11.5 and 4.2]{KS}. Each
corepresentation of $\AA$ can be obtained from  $1$, $u$, $\DD$, and $\dett$
(resp.\  $1$ and $u$ for $\oslqn[2]$) by taking
tensor products, direct sums, and formal differences. 
Using suitable normalisation, $\zeta(v\ot w)=\zeta(v)\zeta(w)$ and
$\zeta(v\oplus w)= \zeta(v)+\zeta(w)$ for irreducible corepresentations $v$ and
$w$. Consequently, $\AA_\inv=\C[U,\DD,\dett]$ (resp.\  $\AA_\inv=\C[U]$ for
$\oslqn[2]$). Next we will show that $\RR':=\{x,\wt{\DD},
\wt{\dett}\}\AA_\inv$ (resp.\  $\RR':=x\AA_\inv$ for $\oslqn[2]$)
coincides with $\RR_\inv$. For we will prove that
\begin{equation}
\AA_\inv=\RR'+\langle U\rangle+\C.
\end{equation}
Since $\DD a=\wt{\DD} a+a\in\RR'+a$ and $\dett a=\wt{\dett}a+a\in\RR'+a$,
$a\in\AA_\inv$, it suffices to prove that $U^k\in\RR'+\langle U\rangle+
\C$,
$k\in\NO$. For $k=0,1$  this is trivial. For $k\ge 2$ use
$U^2=x+(\lam+\fr{s})U-\lam\fr{s}$ and induction on $k$.
Moreover \rf[lin] implies 
$\AA_\inv=\RR_\inv\oplus\langle U\rangle\oplus\C$. 
Hence by the above equation, $\codim_{\AA_{\scriptscriptstyle\inv}}\RR'\le2=
\codim_{\AA_{\scriptscriptstyle\inv}}\RR_\inv$. 
Since $\RR'\subseteq\RR_\inv$, $\RR'=\RR_\inv$. 
\\
By \rf[e-vtvtb] and \rf[e-sub], $\eta\rac\DD^{\pm1}=\eta$ and
$(\vt\ott\vt)\rac\DD^{\pm1}=\vt\ott\vt$. Hence $\xi\rac\DD^{\pm1}=\xi$. 
By the preceding and since $\SS(ra)=\SS(r)\rac a$, $r\in\RR$, $a\in \AA$, we
have $\SS(\RR_\inv)=\SS(\RR')=\SS(x)\rac\AA_\inv=\xi\rac\AA_\inv$. Consequently,
$\SS(\RR_\inv)=\xi\rac\C[U]$ because $\DD$ and $\dett$ act trivial 
on $\xi$. 
Thus $\SS(\RR_\inv)=\langle\xi\rangle$ since $\xi\rac U=(A+\lam)\xi$ by
\rf[xxx]. This completes the proof of Step~1.
\\
{\em Step~2}. All arguments go through for both $\AA=\oglqn[2]$ and
$\AA=\oslqn[2]$. We will show that $\SS$ acts injectively 
on each of the subspaces  spanned by the following 
five resp.\  three elements of   $\RR$:
\begin{align*}
&(u^1_1)^2+q^2(u^2_2)^2-(1+q^2)(u^1_1u^2_2 +\qm u^1_2u^2_1),\, (u^1_2)^2,\,
(u^2_1)^2,\, u^1_2(u^1_1-u^2_2),\, u^2_1(u^1_1-u^2_2),
\\
&(U-z\fr[+]{r})u^1_2, \,(U-z\fr[+]{r})u^2_1,\,(U-z\fr[+]{r})(u^1_1-u^2_2),
\end{align*}
cf.\ \cite[(1.24) and (1.25)]{Wo2}. Since both subspaces are minimal with
respect to $\adr$ and since $\SS$ intertwines $\adr$ and $\dr$, it suffices to
prove that $\SS$ is non-zero on each of these subspaces. By Lemma\,\ref{l-ql}\,(ii) and since  $(P_-)^{11}_{xy}=0$,
  $(P_+)^{11}_{xy}=\delta_{x1}\delta_{y1}$, 
  $(P_+)^{ce}_{22}=\delta_{c2}\delta_{e2}$, and
  $(P_+)^{1z}_{2d}=\QP^{-1}\delta_{z2}\delta_{d1}$ we obtain 
\begin{align*}
\SS(u^1_2u^1_2)&=\sum_{\eps,\iota,\iota'}c_{\eps,\iota}
c_{\eps,\iota'}(P_\eps)^{11}_{xy}(P_\iota)^{mn}_{yz}
(P_\eps)^{xz}_{cd}(P_{\iota'})^{vw}_{dl}
(P_\eps)^{ce}_{22}\om_{mn}\ot\om_{vw},
\\
&=\sum_{\iota,\iota'}c_{+,\iota}
c_{+,\iota'}(P_\iota)^{mn}_{1z}
(P_+)^{1z}_{2d}(P_{\iota'})^{vw}_{d2}\om_{mn}\ot\om_{vw},
\\
&=\sum_{\iota,\iota'}\QP^{-1}c_{+,\iota}
c_{+,\iota'}(P_\iota)^{mn}_{12}
(P_{\iota'})^{vw}_{12}\om_{mn}\ot\om_{vw}=\QP^{-1}c_{+,1}^2
\om_{12}\ot\om_{12}=\QP\om_{12}\ot\om_{12}.
\end{align*}
Since $q^2+1\ne0$, $\SS$ is injective on the $5$-dimensional subspace.
\\
In a  similar way we will compute some part of 
$\SS\bigl((U-\lam)u^1_2\bigr)=q^{-2}\SS(u^1_1u^1_2)
+q^{-4}\SS(u^2_2u^1_2)-\lam(\om_{11}\ot\om_{12}+\om_{12}\ot\om_{22})$.
With the abbreviation $(\om_\iota)_{ij}=\sum_{mn}(P_\iota)_{ij}^{mn}\om_{mn}$
one gets
\begin{align*}
\SS(u^1_1u^1_2)&=\sum_{\iota,\iota'}\QP^{-2}c_{+,\iota}c_{+,\iota'}\bigl(
(\om_\iota)_{12}\ot(\om_{\iota'})_{11}+q^2(\om_\iota)_{12}\ot(\om_{\iota'})_{22}+q\QP(\om_{\iota})_{11}\ot(\om_{\iota'})_{12}\bigr),
\\
\SS(u^2_2u^1_2)&=\sum_{\iota,\iota'}\QP^{-2}c_{+,\iota}c_{+,\iota'}\bigl(
(\om_\iota)_{22}\ot(\om_{\iota'})_{12}+q^{-2}(\om_\iota)_{11}\ot(\om_{\iota'})_{12}+q^{-1}\QP(\om_{\iota})_{12}\ot(\om_{\iota'})_{22}\bigr).
\end{align*}
For the coefficient of $\om_{12}\ot\vt$ in
$\SS\bigl((U-\lam)u^1_2\bigr)$ with respect to the basis
$\{\om_{12},\,\om_{21},\,\om_{11}-\fr{s}^{-1}\vt,\,\vt\}$
we get
\begin{equation*}
q^2\QP^{-1}(q^2\pm q+1)^{-1}(\qm+\qm[5])(q^4\pm q^3+q^2 \pm q +1).
\end{equation*}
Since $q$ is not a root of unity, $\SS\bigl((U-\lam)u^1_2\bigr)$ is
nonzero. Hence $\SS$ is injective on the $3$-dimensional subspace. The proof
of $(\ker\ant[2])_\ll=\SS(\RR)\oplus\langle\vt\ot\vt\rangle$ is complete. To
the last assertion. Obviously $2A-\fr{n}B=0$. Since $\xi=A\eta+B\vt\ot\vt$,
$\eta=-2\fr{n}^{-1}\vt^2$ in $\gdu$. By \rf[e-trans] and by $\eta=-\dd\vt$,
$\dd\theta=2\theta^2$. Thus by Lemma\,\ref{l-inner}\,(ii), $\theta^2$ is central and $\gdu$ is inner; this completes the proof.
\\[1cm]
{\em Acknowledgement.}\/ The author is grateful to I.\,Heckenberger, S.\,Majid,
A.\,Ram and A.\,Sudbery for valuable discussions.

\end{document}